\documentclass{article}
\usepackage{latexsym}
\usepackage{amsmath}
\usepackage{amsfonts}
\usepackage{amssymb}
\usepackage{graphicx}
\newtheorem{theorem}{Theorem}

\newtheorem{corollary}[theorem]{Corollary}

\newtheorem{definition}[theorem]{Definition}

\newtheorem{lemma}[theorem]{Lemma}
\newtheorem{notation}[theorem]{Notation}

\newtheorem{proposition}[theorem]{Proposition}
\newtheorem{remark}[theorem]{Remark}

\begin{document}

\author{Bong H. Lian\\Brandies University\\Department of Mathematics\\Waltham, MA 02254
\and Andrey Todorov\\UC Santa Cruz\\Department of Mathematics\\Santa Cruz, CA 95064\\Institute of Mathematics\\Bulgarian Academy of Sciences\\Sofia, Bulgaria
\and Shing-Tung Yau\\Harvard University\\Department of Mathematics\\Cambridge, MA 02138}
\title{Maximal Unipotent Monodromy for Complete Intersection
CY Manifolds}
\date{March 2000}
\maketitle
\begin{abstract}
The computations that are suggested by String Theory in the B model requires
the existence of degenerations of CY manifolds with maximum unipotent monodromy.
In String Theory such a point in the moduli space is
called a large radius limit (or large complex structure limit). 
In this paper we are going to construct one
parameter families of $n$ dimensional Calabi-Yau manifolds, which are complete
intersections in toric varieties and which have a monodromy operator $T$ such
that (T$^{N}-id)^{n+1}=0$ but (T$^{N}-id)^{n}\neq0,$ i.e the monodromy
operator is maximal unipotent.
\end{abstract}
\tableofcontents

\section{Introduction}

\subsection{General Remarks}

One of the most important problems in algebraic geometry is the study of
families of algebraic varieties parameterized by a variety.
Of special interest is the subvariety in the parameter space
that parameterizes the singular fibers. This subvariety is called the
discriminant locus. One of the main invariants of the discriminant locus is
the so called monodromy group. The monodromy group is defined by
the action of fundamental group of
the complement of the discriminant locus on the cohomology ring of a fixed
non-singular fiber. Of special interest is the action of the monodromy group
on the middle cohomology group. The structure of such actions is of profound
importance. \ In number theory, the counterpart is the action of Galois group
on \`{e}tale cohomology.

In this paper we are going to study the simplest case of the above described
setting, namely we are going to study families of algebraic manifolds over the
unit disk $\mathcal{D}$. We will assume that the only singular fiber is over the center of
the disk, i.e. over $0\in\mathcal{D}.$ From the discussion above it follows
that we obtain a representation of the fundamental group $\pi_{1}%
(\mathcal{D}$
$\backslash$%
$0)=\mathbb{Z}$ in the middle homology, i.e. in $H_{n}(X_{t},\mathbb{Z}%
$)/$Tor$. The finite dimensional representations of $\mathbb{Z}$ are
classified by the Jordan blocks of the linear operator corresponding to
$1\in\mathbb{Z}$. We will give a complete answer to the structure of the
monodromy operator in terms of the topology of the singular fiber. Our method
of proof is based on  Clemens' theory of monodromy and 
the theory of mixed Hodge structures. We find a simple criterium 
for the monodromy operator to have a Jordan block of maximal rank.
This criterium is based on Leray's theory of residues. We will apply this
simple criterium to very concrete examples of complete intersections in
$\mathbb{CP}^{N}$ and toric varieties.

The existence of such degenerations is of prime importance in mirror geometry
and in string theory. The computations that were suggested by String Theory in
the B model required the existence of degenerations of CY manifolds with
maximal unipotent monodromy. In String Theory such a
point in the moduli space is called a large radius (complex structure) limit. The case of
hypersurfaces in toric varieties was treated in \cite{HLY}, where
the construction of a point of maximal degeneracy is done by
studying the GKZ hypergeometric system governing periods of the hypersurfaces.

\ \ For recent important developments in mirror geometry see \cite{LLY1},
\cite{LT} and references therein. The results of this paper are closely
related to Strominger-Yau-Zaslow conjecture.

\subsection{Description of the Paper}

In \textbf{Section 2} we introduce the basic notions and review some results
stated in \cite{KKMS}.

In \textbf{Section 3} we review the generalization of Lefschetz's theory of
vanishing cycles due to Clemens. We describe Clemens' method for constructing
Jordan blocks in the monodromy operator.

In \textbf{Section 4} we prove a general formula for the number of Jordan
blocks in the monodromy operator in terms of some invariants of the singular
fiber of one parameter family of K\"{a}hler manifolds.

In \textbf{Section 5} we prove a simple criterium for the existence of a
Jordan block of maximum size  in the monodromy operator in terms of
Leray's residue calculus.
We also construct families of complete intersections of CY
manifolds in $\mathbb{CP}^{N}$ whose monodromy operators contain Jordan block
of maximum size.

In \textbf{Section 6} we review some basic facts in toric geometry and
construct families of complete intersections of CY manifolds in toric
varieties whose monodromy operators contain Jordan block of maximum size
generalizing the construction in section 5.

In \textbf{Section 7} we briefly discuss the connection between the
present approach and a previous approach which uses hypergeometric functions.
We also discuss an interesting relationship between maximal unipotent monodromy
and the SYZ conjecture, and illustrate this in the case of polarized K3 surfaces.

In the Appendix, we give a complete description of Clemen's cell complex
for hypersurfaces in a toric variety.

{\it Acknowledgements:}
B.H.L.'s research
is supported by NSF grant DMS-0072158.
S.T.Y.'s research is supported by DOE grant
DE-FG02-88ER25065 and NSF grant DMS-9803347.

\section{Basic Definitions and Notations}

\subsection{Mumford's Semi-Stable Reduction Theorem}

In this article we study one-parameter families of $n$ dimensional K\"{a}hler manifolds

\begin{center}
$\pi:\mathcal{X}\rightarrow\mathcal{D}$
\end{center}

over a disk $\mathcal{D}$. We assume that $\mathcal{X}$ is a smooth algebraic
manifold and that for each $t\neq0$ $\pi^{-1}(t)=X_{t}$ is a non-singular
$n$ dimensional K\"{a}hler manifold. From Hironaka's theorem on the
resolution of singularities, we may assume that the singular fibre $X_{0}%
=\pi^{-1}(0)$ is a divisor with normal crossings $X_{0}=\cup C_{i}. $

We will choose local coordinates in the following manner. Let $x_{n+1}\in
C_{i_{1}}\cap...\cap C_{i_{n+1}}\subset X_{0}.$ Let $U_{x_{n+1}}$ be an open
polycylinder in $\mathcal{X}$ containing the point $x_{n+1}.$ Let 
$z_{i_j}=0$ be the defining equation of the divisor $C_{i_{j}}$ in
$U_{x_{n+1}}.$ It is easy to see that after shrinking the disk
$\mathcal{D}$, the fibers of the map $\pi:\mathcal{X}\cap U_{x_{n+1}}%
\rightarrow\mathcal{D}$ are locally given by

\begin{center}
$z_{i_{1}}^{m_{1}}...z_{i_{n+1}}^{m_{n+1}}=t,$ $m_{j}\geq1$ $\&$ $m_{j}%
\in\mathbb{Z}$.
\end{center}

In the same manner, let $x_{k}\in X_0$ and
let $C_{i_{1}},..,C_{i_{k}}$ be those components of $X_0$ containing $x_k$.
Let $U_{x_{k}}$ be an open polycylinder in $\mathcal{X}$
containing $x_k$ but not intersecting $C_j$ with $j\neq i_1,..,i_k$.
Then the fibers of
the map $\pi:\mathcal{X}\cap U_{x_{k}}\rightarrow\mathcal{D}$ are locally given by

\begin{center}
$z_{i_{1}}^{m_{1}}...z_{i_{k}}^{m_{k}}=t,$ $m_{j}\geq1$ $\&$ $m_{j}%
\in\mathbb{Z}$.
\end{center}


In \cite{KKMS} Mumford proved\ that after taking a finite covering of the disk
$\mathcal{D}$, lifting the family and resolving the singularities, we may
assume that the the fibres of the map $\pi$:$\mathcal{X}\rightarrow
\mathcal{D}$ are given locally by $z_{1}^{k_{1}}...z_{n+1}^{k_{n+1}}=t,$
where $k_{j}$ is either 0 or 1. From now on we are going to
assume that we are in the above setting justified by Mumford's theorem.

\subsection{Geometric and Homological Monodromy (Basic Properties)}

If we restrict our family $\pi:\mathcal{X}\rightarrow\mathcal{D}$ to the
circle $S^{1}:=\partial\mathcal{D}$ then we get a representation of $\pi
_{1}(S^{1})=\mathbb{Z}$ \ in the group of diffeomorphisms of $X_{t}.$ Indeed
if we remove a point from $S^{1}$,and restrict our family to $S^{1}$
$\backslash$ $s$, $\pi_{1}:\mathcal{X}|_{S^{1}}\rightarrow S^{1}$ where $s\in
S^{1},$ we will get a trivial $C^{\infty}$ family $\left(  S^{1}\backslash
s\right)  \times X_{t}\rightarrow S^{1}$ $\backslash$ $s.$ Our family $\pi
_{1}:\mathcal{X}|_{S^{1}}\rightarrow S^{1}$ is obtained from the trivial
family $\left(  S^{1}\backslash s\right)  \times X_{t}\rightarrow S^{1}$
$\backslash$ $s$ by ''gluing'' it by $\phi\in Diff(X_{t}).$ The diffeomorphism
$\phi$ is the generator of \ $\pi_{1}(S^{1})=\mathbb{Z}.$ We will call $\phi$
the geometric monodromy. The induced action of $\phi$ on $H_{n}(X_{t}%
,\mathbb{Z})$ will be called the monodromy operator and will be denoted by
$T.$ The main result about the operator $T$ is that we have always:

\begin{center}
$\left(  T^{N}-id\right)  ^{n+1}=0$
\end{center}

for some positive integer $N$. Here $n$ is the complex dimension of $X_{t}.$
This theorem was proved by many mathematicians including Griffiths, N. Katz,
Clemens, Landesman, Deligne and so on.

From Mumford's result, we can assume that the following conditions
hold for the family $\pi:\mathcal{X}\rightarrow\mathcal{D}:$

\textbf{1.} $\pi^{-1}(0)=X_{0}=\cup_{i=0}^{m}C_{i}$ is a divisor of normal
crossings and for $k=0,..,n,$ $C_{i_{0}}\cap..\cap C_{i_{k}}%
,i_{0}<...<i_{k},$ is a non-singular irreducible subvariety of Sing$_{n-k}%
$($X_{0})$ if non-empty.

\textbf{2.} The fibers of the map $\pi:\mathcal{X}\rightarrow\mathcal{D}$ are
locally given in the open policylinders $\{U\}$ defined above by \ $z_{i_{1}%
}^{k_{1}}...z_{i_{n+1}}^{k_{n+1}}=t,$ where $k_{j}$ is either 0 or 1.

\begin{definition}
\label{Mum}If the family $\pi:\mathcal{X}\rightarrow\mathcal{D}$ satisfies the
above conditions we will say that it is in normal form.
\end{definition}

From now on we will consider only families $\pi:\mathcal{X}\rightarrow
\mathcal{D}$ in normal form.

\begin{notation}
\label{no}Given $X_{0}=C_{0}\cup...\cup C_{m}$ we shall use the following notations:
\end{notation}

\begin{center}
$I=\{i_{0},...i_{k}\}$ is an index set with $i_{0}<...<i_{k}$

$C_{I}=C_{i_{0}}\cap...\cap C_{i_{k}},$ $|I|=k+1;$ and

$C^{[k]}=\underset{|I|=k+1}{\cup}C_{I}.$
\end{center}

\section{Review of Clemens' Theory of Geometric Monodromy}

\subsection{Construction of Clemens' retraction map}

Let $\pi:\mathcal{X}\rightarrow\mathcal{D}$ \ be a family of K\"{a}hler
manifolds as defined in Definition \ref{Mum}. We will construct a contraction map:

\begin{center}
$\mathcal{C}:\mathcal{X}\rightarrow X_{0}.$
\end{center}

\subsubsection{Construction of Clemens' Vector Field}

The local description of our family \ $\pi:\mathcal{X}\rightarrow\mathcal{D}$ 
is given in $\mathbb{C}^{n+1}$ by the equation:

\begin{center}
$z_{1}...z_{k}=t$
\end{center}
For some $k=1,..,n+1$.
Without loss of generality we may assume that $t$ is a real number.

Let $z_{i}=r_{i}e^{2\pi\sqrt{-1}\phi_{i}},$then the equation $z_{1}...z_{k}=t$
for $k=1,..,n+1$ is equivalent to the equations:

\begin{center}
$r_{1}...r_{k}=t$ and $\sum_{i=1}^{k}\phi_{i}=0$ for $k=1,..,n+1.$
\end{center}

We will construct first a local vector field which will define the retraction
map in this local situation and then using partition of unity we will
construct a global vector field and thus the retraction map. It is easy to see
that it is enough to construct a retraction of hyperbola in $\ \mathbb{R}%
_{+}^{n+1}\ $\ given by

\begin{center}
$\ r_{1}...r_{n+1}=t$
\end{center}

to the union of coordinate hyperplanes $r_{1}=0,...,r_{n+1}=0,$ where
$r_{i}>0$ and $t>0.$ Let us suppose that we will consider only hyperbolas
$r_{1}...r_{n+1}=t$ for $0<t<\frac{1}{2}.$ We will construct special vector
field in $\ \mathbb{R}_{+}^{n+1}.$ We will need to define some notions.

We will denote by $\mathcal{Q}$ the following set in $\mathbb{R}_{+}^{n+1}$

\begin{center}
$\mathcal{Q}$:=$\left\{  r=\left(  r_{1},..,r_{n+1}\right)  |\,r_{i}%
\geq0\text{ and\ there exists }1\leq i\leq n+1\text{ such that 0%
$<$%
r}_{i}<1\right\}  .$
\end{center}

$\mathcal{Q}(1)$ will be the unit cube in $\mathbb{R}_{+}^{n+1},$ i.e.

\begin{center}
$\mathcal{Q}(1):=\left\{  r=\left(  r_{1},..,r_{n+1}\right)  |0\leq r_{i}%
\leq1\text{ for all }1\leq i\leq n+1\right\}  $
\end{center}

It is clear that for $0<t<\frac{1}{2}$ the hyperbolas $r_{1}...r_{n+1}=t$ are
contained in $\mathcal{Q}$.

\subsubsection{Construction of the Vector Field in \ $\mathbb{R}_{+}^{n+1}$}

Let us suppose that $r=\left(  r_{1},..,r_{n+1}\right)  \in\mathcal{Q}$ is a
point such that $0\leq r_{i_{1}}\leq1,...,0\leq r_{i_{k}}\leq1.$ To the point
$r\in\mathcal{Q}$ $\ $we will assign a point $r(i_{1},..,i_{k})=(r_{1}%
,...,r_{n+1})\in\mathbb{R}_{+}^{n+1}$ $r_{i_{1}}=...=r_{i_{k}}=1$ and the rest
of the coordinates of \ $r(i_{1},..,i_{k})$ are the same as the point
$r\in\mathcal{Q}.$ Let $l(r)$ be the line that joints the point $r\in
\mathcal{Q\subset}\mathbb{R}_{+}^{n+1}r(i_{1},..,i_{k})\in\mathbb{R}_{+}%
^{n+1}$ with the point $r(i_{1},..,i_{k})\in\mathbb{R}_{+}^{n+1}.$ In this way
we define a vector field in $\mathcal{Q\subset}\mathbb{R}_{+}^{n+1}.$ Using a
partition of unity and using the vector field in $\mathcal{Q}$ that we have just
defined, we obtain a global vector field on $\pi:\mathcal{X\rightarrow D}.$

When we integrate this vector field, we obtain the Clemens' contraction map
$\mathcal{C}:\mathcal{X\rightarrow}X_{0}.$ Thus we obtain a map for each
$t\in\mathcal{D}$ $\backslash$ $0:\mathcal{C}_{t}:X_{t}\rightarrow X_{0}.$

\begin{definition}
\label{nik}The map $\mathcal{C}_{t}$ as defined above will be called Clemens map.
\end{definition}

For details of the construction see \cite{Cl}.

\subsection{Properties of the Clemens' map}

It is easy to prove the $\mathcal{C}_{t}$ has the following properties:

\begin{lemma}
\label{CL0} The Clemens' map $\mathcal{C}_{t}$ has the following properties:
\textbf{(i)} Suppose $z\in C_{I}$, then $\mathcal{C}_{t}^{-1}(z)=\left(
S^{1}\right)  ^{k}$ is a k dimensional real torus and \textbf{(ii)
}$\mathcal{C}_{t}$ defines a diffeomorphism between $X_{t}$
$\backslash$%
$\mathcal{C}_{t}^{-1}(\operatorname{Sing}(X_{0}))$ and $X_{0}$
$\backslash$%
Sing($X_{0}).$
\end{lemma}

\textbf{Proof\ of\ Lemma\ }\ref{CL0}: For the proof of Lemma \ref{CL0} see
\cite{Cl}. $\blacksquare.$

\begin{definition}
\label{t1} Let $\mathcal{T}(i_{1},..,i_{k})$ be the tubular neighborhood of
$C_{i_{0}}\cap..\cap C_{i_{k}}$ in $C_{i_{1}}\cap..\cap C_{i_{k}}.$ We will
denote by $p(i_{1},..,i_{k}):\mathcal{T}(i_{1},..,i_{k})\rightarrow C_{i_{0}}\cap
C_{i_{1}}\cap...\cap C_{i_{k}}$ the projection maps for any $i_{1},..,i_{k}$.
\end{definition}

\begin{definition}
\label{t10}Given an $(n-k)$-cycle $\gamma\in H_{n-k}(C_{i_{0}}\cap...\cap C_{i_{k}%
},\mathbb{Q})$, we define an $n$-cycle $p_{k}^{-1}(\gamma)\in H_{n}%
(X_{0}\backslash$Sing($X_{0}))$ as follows.
Let 
$$p_{1}
^{-1}(\gamma):=\partial(p(i_{1},..,i_{k})^{-1}(\gamma)).$$
This is
the boundary of $p(i_1,..,i_k)^{-1}(\gamma)$, hence it is 
cycle of dimension $n-k+1$ in
$C_{i_{1}}\cap...\cap C_{i_{k}}\backslash C_{i_{0}
}\cap...\cap C_{i_{k}}$. 
Let
$$p_{2}^{-1}(\gamma)=\partial(p(i_{2},..,i_{k})^{-1}p_1^{-1}(\gamma)).$$
This is a cycle of dimension $n-k+2$ in
$C_{i_{2}}\cap...\cap C_{i_{k}}\backslash C_{i_{1}
}\cap...\cap C_{i_{k}}$. 
By continuing this way,
we define at the end 
$$p_{k}^{-1}(\gamma)=\partial(p(i_{k})^{-1}p_{k-1}^{-1}(\gamma)).$$
This is a cycle of dimension $n$ in
$C_{i_{k}}\backslash C_{i_{k-1}}
\cap C_{i_{k}}$.
We denote by $\pi_{k}:p_{k}^{-1}(\gamma)\rightarrow
\gamma$ the natural projection.
\end{definition}

\begin{proposition}
\label{t11}For any point $z\in\gamma\in H_{n-k}(C_{i_{0}}\cap..\cap C_{i_{k}%
},\mathbb{Q})$ we have $\pi_{k}^{-1}(z)=\left(  S^{1}\right)  ^{k},$ where
$\pi_{k}:p_{k}^{-1}(\gamma)\rightarrow\gamma$ is defined above. In other words
the preimage of a point is a $k$ dimensional real torus.
\end{proposition}

\textbf{Proof\ of\ Proposition }\ref{t11}\textbf{:} Proposition \textbf{\ }%
\ref{t11} follows directly from Definition \ref{t10} of the cycle
$p_{k}^{-1}(\gamma)\in H_{n}(X_{0}\backslash$Sing($X_{0}))$ . $\blacksquare.$

\begin{definition}
\label{t2}Given $\gamma\in H_{n-k}(C_{i_{0}}\cap..\cap
C_{i_{k}},\mathbb{Q})$, we denote $\gamma_{t}:=\mathcal{C}_{t}%
^{-1}(\gamma)$ and $\gamma_{t}(j):=\mathcal{C}_{t}^{-1}\left(
p_{j}^{-1}\left(  \gamma\right)  \right)$, which are $n$-cycles
representing elements
in $H_{n}(X_{t}
,\mathbb{Q})$.
\end{definition}

\begin{lemma}
\label{an1} The cycles $\gamma_{t}$ and $\gamma_{t}(j)$ above
are homological to each other for 1$\leq j\leq k$.
\end{lemma}

\textbf{Proof\ of\ Lemma \ref{an1}:} We will prove first that $\gamma_{t}(j) $
and $\gamma_{t}(j+1)$ are homological to each other for $1\leq j\leq k$,
i.e. there exists a $n+1$ chains $\Gamma_{0}(j)$ in $C_{i_{j+1}}\cap..\cap
C_{i_{k}}\backslash C_{i_{j}}\cap..\cap C_{i_{k}}\subset X_{0}$ as follows
using Definitions \ref{t1} and \ref{t10}: $\Gamma_{0}(j):=p(i_{j}%
,..,i_{k})^{-1}(p_{j}^{-1}(\gamma))\backslash(p_{j}^{-1}(\gamma)).$ \ Let us
define the $n+1$ dimensional chain $\Gamma_{t,\gamma}\left(  j\right)  $ in
$X_{t}$ as follows: $\Gamma_{t,\gamma}\left(  j\right)  :=\mathcal{C}_{t}%
^{-1}(\Gamma_{0}(j)).$ Lemma \ref{CL0} implies directly that 

\begin{center}
$\partial\left(  \Gamma_{t,\gamma}\left(  j\right)  \right)  =\mathcal{C}%
_{t}^{-1}(p_{j}^{-1}(\gamma))-\mathcal{C}_{t}^{-1}(p_{j-1}^{-1}(\gamma
))=\gamma_{t}(j)-\gamma_{t}(j-1).$
\end{center}

\ Lemma \ref{an1} is proved. $\blacksquare.$

In \cite{Cl} the following result was proved:

\begin{theorem}
\label{cl1}Let $\gamma\in H_{n-k}($Sing$_{k}(X_{0}),\mathbb{Z})$ \textit{be
such that }$\mathcal{C}_{t}^{-1}(\gamma)=\gamma_{t}\in H_n(X_t,\mathbb{Z})$ 
\textit{is a non-zero.
Then there exists cycles }
$\alpha_1,..,\alpha_k\in H_n(X_t,\mathbb{Z})$ \textit{such that }$T(\alpha_j)
=\gamma_t+\sum_{i=1}^{j}\alpha_i$ \textit{for 1}$\leq j\leq k$.
\end{theorem}
We sketch the construction here.
(See \cite{Cl} for details.)

{\it Clemens' construction of Jordan block by Picard-Lefschetz
Duality.} Let $\gamma\in H_{n-k}(C_{i_{0}}\cap...\cap C_{i_{k}},\mathbb{Z})$ be
a cycle such that $\mathcal{C}_{t}^{-1}(\gamma)$ be a non-zero element in
$H_{n}(X_{t},\mathbb{Z\dot{)}}.$ In Definition \ref{t10} we defined a cycle
$p_{1}^{-1}(\gamma)$ in $C_{i_{2}}\cap...\cap C_{i_{k}}\backslash C_{i_{1}%
}\cap...\cap C_{i_{k}}.$ It is easy to see by using the fact that
$\mathcal{C}_{t}^{-1}(\gamma)\in H_{n}(X_{t}%
,\mathbb{Z\dot{)}}$ is nonzero and Lemma \ref{an1} that

\begin{center}
$p_{1}^{-1}(\gamma)\in$ $H_{n-k+1}(C_{i_{2}}\cap...\cap C_{i_{k}}\backslash
C_{i_{1}}\cap...\cap C_{i_{k}},\mathbb{Z})~\&~p_{1}^{-1}(\gamma)\neq0.$
\end{center}

Let us denote by $\gamma_{1}\in H_{n-k+1}(C_{i_{2}}\cap...\cap C_{i_{k}%
};C_{i_{1}}\cap...\cap C_{i_{k}},\mathbb{Z})$ the Picard Lefschetz dual cycle
of $p_{1}^{-1}(\gamma).$ Let $T_{i_{3},...,i_{k}}(\overline{\gamma_{1}})$ be
the tubular neighborhood of the closure of $\overline{\gamma_{1}}$ in
$C_{i_{3}}\cap...\cap C_{i_{k}}\backslash C_{i_{2}}\cap...\cap C_{i_{k}}.$ Let
us denote by $p_{2}^{-1}(\overline{\gamma_{1}})$ the boundary of
$T_{i_{3},...,i_{k}}(\overline{\gamma_{1}}),$ i.e.

\begin{center}
$p_{2}^{-1}(\overline{\gamma_{1}})$ $=\partial T_{i_{3},...,i_{k}}%
(\overline{\gamma_{1}})$
\end{center}

It is easy to see that

\begin{center}
$p_{2}^{-1}(\overline{\gamma_{1}})\in H_{n-k+2}(C_{i_{3},...,i_{k}}\backslash
C_{i_{2},...,i_{k}},\mathbb{Z}).$
\end{center}

Let us denote by $\gamma_{3}\in H_{n-k+2}(C_{i_{3},...,i_{k}};C_{i_{2}%
,...,i_{k}},\mathbb{Z})$ the Picard-Lefschetz dual to $p_{2}^{-1}%
(\overline{\gamma_{1}}).$ We can continue this process and thus we will define
cycles $\gamma,\gamma_{1},...,\gamma_{k},$ where $\gamma_{j}\in H_{n-k+2}%
(C_{i_{j},...,i_{k}};C_{i_{j-1},...,i_{k}},\mathbb{Z}).$ Clemens proved in
\cite{Cl}\ that the monodromy operator \textit{T} acts as follows on
$\mathcal{C}_{t}^{-1}(\gamma),
\mathcal{C}_{t}^{-1}(\gamma_1),...,
\mathcal{C}_{t}^{-1}(\gamma_{k}):$

\begin{center}
T($\mathcal{C}_{t}^{-1}(\gamma))=\mathcal{C}_{t}^{-1}(\gamma),...,$%
T($\mathcal{C}_{t}^{-1}(\gamma_{k}))=\mathcal{C}_{t}^{-1}(\gamma)+\sum
_{j=1}^{k}\mathcal{C}_{t}^{-1}(\gamma_{j}).$
\end{center}

\begin{corollary}
\label{betty2} Let $\pi:\mathcal{X}\rightarrow\mathcal{D}$ be a family of
K\"{a}hler manifolds over the disk such that:
\end{corollary}

\begin{enumerate}
\item \textit{For t}$\neq0,$ $\pi^{-1}(t):=X_{t}$ \textit{is a non singular
variety of complex dimension }$n\geq1.$

\item $\pi^{-1}(0)=X_{0}=\cup_{i=0}^{m}C_{i}$ \textit{is a divisor of normal
crossing and } $\pi$ \textit{is locally given by} $z_{1}^{n_{1}}%
\times...\times z_{k}^{n_{k}}=t,$ \textit{where} $n_{i}$ \textit{are positive integers.}

\item \textit{\ Suppose }$C_{0}\cap..\cap C_{n}$\textit{\ is a point}
\textit{and }$\mathcal{C}_{t}^{-1}(C_{0}\cap..\cap C_{n})=\gamma_{t}$
\textit{is a non zero cycle in} $H_{n}(X_{t},\mathbb{Q}),$
\end{enumerate}

\textit{then the monodromy operator of the family }$\pi:\mathcal{X}%
\rightarrow\mathcal{D}$\textit{\ has a Jordan block of size }$n+1.$
(See \cite{Cl}.)

\section{The Jordan Normal Form of the Monodromy Operator}

We begin by introducing the following combinatorial invariant of a family of
algebraic varieties $\pi:\mathcal{X}\rightarrow\mathcal{D}$ put in a Mumford form.
This will be needed later.

\begin{definition}
\label{POL} We will define Clemens' simplicial complex of the family
$\pi:\mathcal{X}\rightarrow\mathcal{D}$ as follows: To each divisor $C_{i}$ we
attached a point $p_{i}$ in $\mathbb{R}^{d}$, where d is a large enough
integer. We will assume that the points p$_{i}$ are in general position, i.e.
they do not lie in a hyperplane. If $C_{i}$ intersects $C_{j}$ then we attach
to the points p$_{i}$ and p$_{j}$ one dimensional simplex. If $C_{i}$, $C_{j}$
and $C_{k}$ intersect then we attached a two dimensional simplex on p$_{i}$,
p$_{j}$ and p$_{k}.$ We continue in that manner and we obtain a simplicial
complex $\Pi(X_{0})$ which we will call Clemens' simplicial complex.
\end{definition}

\subsection{Definition of the Gysin Map}

\begin{definition}
\label{G}
Let $X$ be a compact complex manifold, and $C$ be a divisor of normal crossing.
The map 
$G_k:H_{k+2}(X,\mathbb{Z})\rightarrow H_{k}(C,\mathbb{Z})$ defined by
$G_{k}(\gamma):=\gamma\cap\lbrack C],$ where $\gamma\cap\lbrack C]$ means
intersection of class of cohomology in $X$, is called the Gysin map.
\end{definition}

\begin{remark}
We will use the Gysin map in case $C_{i_{0}}\cap...\cap
C_{i_{k}}\subset C_{i_{1}}\cap..\cap C_{i_{k}}$ and will denote by
$$
G_{k}:\underset{i_{1},..,i_{k}}{\oplus}H_{n-k+2}(C_{i_{1}}\cap..\cap
C_{i_{k}},\mathbb{Q})\rightarrow H_{n-k}(C_{i_{0}}\cap...\cap C_{i_{k}
},\mathbb{Q)}
$$
where $G_{k}(\gamma)$\textit{\ is the image of the cycle }$\left[
\gamma\cap\left[  C_{i_{0}}\cap..\cap C_{i_{k}}\right]  \right]  $\textit{\ in
}$H_{n-k}($Sing$_{n-1}(X_{0}),\mathbb{Q))}$.
\end{remark}

\begin{definition}
\label{G01}The dual $G_{k}^{\ast}$ of the Gysin map using Poincare duality is
defined as follows for the pair ($X,C)$
\end{definition}

\begin{center}
$G_{k}^{\ast}:H^{k}(C,\mathbb{Z})\rightarrow H^{k+2}(X,\mathbb{Z}),$
\end{center}

\textit{where} $G^{\ast}_k(\alpha)=\alpha\wedge c_{1}[C]$ \textit{and }
$c_{1}[C]$ \textit{is the first Chern class of the line bundle defined by the
normal crossing divisor $C$}.

\subsection{Review of Deligne's Theory of Mixed Hodge Structures}

The cohomology of $X_{0}$
$\backslash$%
Sing($X_{0})$ can be computed as the cohomology of the de Rham log complex
$\mathcal{A}^{\ast}(X_{0},\log<$Sing(X$_{0})>)$. 

\begin{definition}
\label{d1}We will say that a form $\omega$ on one of the components $C_{i}$ of
$X_{0}$ had a logarithmic singularities if for each point $z\in X_{0}$ and
some open neighborhood $U\subset C_{i}$ of the point $z$ we have
\end{definition}

\begin{center}
$\omega|_{U}=\alpha\frac{dz^{1}}{z^{1}}\wedge...\wedge\frac{dz^{k}}{z^{k}},$
\end{center}

\textit{where } $\alpha$ \textit{is a }C$^{\infty}\,$\ \textit{form in }$U$
\textit{e \ and locally on }$X_{0}$ \textit{is defined by the equations }
$z^{1}\times...\times z^{k}=0.$

\begin{definition}
\label{d2}We define the de Rham log complex
as follows:
\end{definition}

\begin{center}
$\mathcal{A}^{\ast}(X_{0},\log<$Sing(X$_{0})>)$

$=\{\omega\in C^{\infty}\left(  X_{0}\backslash\text{Sing(X}_{0}),\Omega
^{\ast}\right)  |\omega$ \textit{and }$d\omega$ \textit{are} C$^{\infty}$
\textit{forms on }$X_{0}\backslash$Sing(X$_{0})$ \textit{which have
log\ singularities on} Sing(X$_{0})\}.$
\end{center}

\begin{remark}
\label{r0} It is easy to see that if $\omega\in\mathcal{A}^{m}(X_{0},\log<
$Sing(X$_{0})>)$ then $\omega$ locally around a point $z\in U\subset C_{i} $
in each of the components of $X_{0}$ and $z\in$Sing($X_{0})$ is given by
\end{remark}

\begin{center}
$\omega|_{U}=\omega_{1}\wedge\frac{dz_{i_{1}}}{z_{i_{1}}}\wedge...\wedge
\frac{dz_{i_{k}}}{z_{i_{k}}},$
\end{center}

\textit{where $\omega_1$ is a C}$^{\infty}$\textit{\ }a $(m-k)$\textit{\ \ form on
}$U\subset C_{i}$;\textit{\ and Sing(}$X_{0})\cap U$\textit{\ \ is given by
}$z_{i_{1}}\cdots z_{i_{k}}=0$\textit{\ in }$U.$\textit{\ }

Deligne proved that there exists a mixed Hodge structure on $X_{0}\backslash
$Sing($X_{0}).$ The existence of Mixed Hodge Structure is based on the
following filtration:

\begin{definition}
\label{Fil} On the complex $\mathcal{A}^{\ast}(X_{0},\log<$Sing(X$_{0})>)$ we
define the weight filtration $W_{l}$ to be those forms $\phi$ that locally
around Sing ($X_{0})$
\end{definition}

\begin{center}
$\phi\in\mathcal{A}^{\ast}(U)\left\{  \frac{dz_{i_{1}}}{z_{i_{1}}}%
\wedge...\wedge\frac{dz_{i_{k}}}{i_{k}}\right\}  .$
\end{center}

\begin{definition}
\label{PR}The Poincare Residue Operator $R^{[k]}:W_{k}\rightarrow
\mathcal{A}^{\ast-k}(C^{[k]})$ is defined by
\end{definition}

\begin{center}
$R^{[k]}\left(  \alpha\wedge\frac{dz_{i_{1}}}{z_{i_{1}}}\wedge...\wedge
\frac{dz_{i_{k}}}{i_{k}}\right)  =\alpha|_{C_{I}}.$
\end{center}

\begin{definition}
\label{ss}Let us consider the decreasing filtration$...\supset W^{-l}\supset
W^{-l+1}\supset...$where $W^{-l}=W_{l}.$ Accordingly there is a spectral
sequence $\{E_{r}\}$ such that $E_{\infty}$ is the associated graded to the
weight filtration in $H^{\ast}(X_{0}\backslash$Sing($X_{0}),\mathbb{C}).$
\end{definition}

The filtration was reversed by Deligne so that we can form a spectral sequence
of the filtered de Rham logarithmic complex. By using the Poincare residue map
Deligne proved the following Theorems:

\begin{theorem}
\label{D0}The cohomology of $X_{0}\backslash$Sing($X_{0})$ are equal to the
cohomology of the De Rham log complex $\mathcal{A}^{\ast}(X_{0},\log
<$Sing(X$_{0})>).$
\end{theorem}

(For the proof of Theorem \ref{D0} see \cite{GS}.)

\begin{theorem}
\label{D}\textbf{\ i. }The spectral sequence defined as above degenerates at
the second step. \textbf{ii. }$E_{1}\backsimeq\oplus H^{\ast}(C_{I})$ and the
mapping $d_{1}:E_{1}\rightarrow E_{1}$ is a morphism of Hodge structures given
by the Gysin map (See Definition \ref{G}.)
\end{theorem}

\begin{center}
$H^{\ast}(C_{i_{0}}\cap...\cap C_{i_{l}})\rightarrow H^{\ast}(C_{i_{1}}%
\cap...\cap C_{i_{l}}).$
\end{center}

\subsection{Jordan Normal Form of the Monodromy Operator}

In this section we will prove the following Theorem:

\begin{theorem}
\label{betty1}\textbf{\ i. }The number of Jordan blocks of rank $k\leq n$ is
equal to the rank of the group
\end{theorem}

\begin{center}
$H_{n-k}(\operatorname{Si}ng_{k}(X_{0}),\mathbb{Q}$)/$\operatorname{Im}\left(
G_{k}\right)  ,$
\end{center}

\textit{where $G_k$ is the Gysin map.}

\textbf{ii.}\textit{The number of Jordan blocks of rank }$n+1$\textit{\ is
equal to} $\dim H_{n}(\Pi(X_{0}),\mathbb{Q}),$ \textit{where }$\Pi(X_{0}%
)$\textit{\ is the Clemens' polyhedra defined in Definition \ref{POL}.}

\textbf{Proof\ of\ Theorem\ \ref{betty1}:}The proof of both part is based on
Corollary \ref{betty2}. Part \textbf{i} of Theorem \ref{betty1} follows
directly from the following three Lemmas and Theorem \ref{cl1}:

\begin{lemma}
\label{G0}\textbf{. }Let $\gamma\in\operatorname{Im}G_{k}\subseteq H_{n-k}%
($Sing$_{k}(X_{0}),\mathbb{Q}$)$,$ then $\gamma_{t}=\mathcal{C}_{t}%
^{-1}(\gamma)$ is homological to zero in $X_{t}.$
\end{lemma}

\begin{lemma}
\label{G1} Suppose that $\gamma\in H_{n-k}($Sing$_{k}(X_{0}),\mathbb{Q}%
$)/$\operatorname{Im}G_{k}$ and $\gamma\notin\operatorname{Im}G_{k\text{ }},$
then there exists a non zero class of cohomology $\widetilde{\omega}\in
H^{n}(X_{0}\backslash$Sing($X_{0}$)) such that 
\end{lemma}

\begin{center}
$\int_{p_{k}^{-1}(\gamma)}\widetilde{\omega}\neq0.$
\end{center}

\textit{where the cycle $p_{k}%
^{-1}(\gamma)\in H_{n}(X_{0}\backslash$Sing($X_{0}$),$\mathbb{Q})$ 
is defined in Definition \ref{t10}.} 

\begin{lemma}
\label{G11} Let $\gamma\in H_{n-k}($Sing$_{k}(X_{0}),\mathbb{Q}$%
)/$\operatorname{Im}G_{k}$ and $\gamma\notin\operatorname{Im}G_{k\text{ }},$
then $\gamma_{t}$ is a non zero element in $H_{n}(X_{t},\mathbb{Q}).$
\end{lemma}

\textbf{Proof\ of\ Lemma \ref{G0}:}We will prove part \textbf{i.} Suppose that

\begin{center}
$\gamma\in\operatorname{Im}G_{k}\subseteq H_{n-k}($Sing$_{k}(X_{0}%
),\mathbb{Q}$).
\end{center}

From the Definition \ref{G} of the Gysin map it follows that there exists a cycle

\begin{center}
$\Gamma\in\underset{i_{1},..,i_{k}}{\oplus}H_{n-k+2}(C_{i_{1}}\cap..\cap
C_{i_{k}},\mathbb{Q})$
\end{center}

such that $\Gamma\cap\underset{i_{0},..,i_{k}}{\oplus}[C_{i_{0}}\cap..\cap
C_{i_{k}}]=\gamma.$ The definition \ref{nik} of the Clemens map it follows
that the boundary of $\mathcal{C}_{t}^{-1}(\Gamma\backslash\gamma) $ is
exactly $\mathcal{C}_{t}^{-1}(\gamma).$ Lemma \ref{G0} is proved.
$\blacksquare.$

\textbf{Proof\ of\ Lemma\ }\ref{G1}: The construction of the form $\omega$ is
based on the following fact about the cohomology (homology) of ($X_{0}$
$\backslash$%
Sing($X_{0})),$ where $X_{0}$ is a K\"{a}hler variety and Sing($X_{0})$ is a
divisor with normal crossings in $X_{0}$.

Suppose that $\gamma\in H_{n-k}($Sing$_{k}(X_{0}),\mathbb{Q}$) and
$\gamma\notin\operatorname{Im}G_{k\text{ }}.$ Let $\gamma\in H_{n-k}(C_{i_{0}%
}\cap...\cap C_{i_{k}},\mathbb{Q}$). In order to construct the form $\omega$
we will need to recall the how the dual $G_{k}^{\ast}$to the Gysin map
$G_{k}:H_{k+2}(X,\mathbb{Z})\rightarrow H_{k}(C,\mathbb{Z})$ is defined in
Definition \ref{G01} as follows:

\begin{center}
$G_{k}^{\ast}:H^{k}(C,\mathbb{Z})\rightarrow H^{k+2}(X,\mathbb{Z}),$ where
$G^{\ast}_k(\alpha)=\alpha\wedge c_{1}[C]$.
\end{center}

$c_{1}[C]$ is the Chern class of the line bundle defined by the divisor with
normal crossings $C$ in X.

We will need the following Proposition:

\begin{proposition}
\label{G101} Let $\gamma\in H_{n-k}(C_{i_{0}}\cap...\cap C_{i_{k}},\mathbb{Q})
$ and $\gamma\notin\operatorname{Im}G_{k}.$ Let $\omega_{n-k}\in
H^{n-k}(C_{i_{0}}\cap...\cap C_{i_{k}},\mathbb{Q}$) and
\end{proposition}

\begin{center}
$\int_{\gamma}\omega_{n-k}\neq0$
\end{center}

\textit{then} $\omega_{n-k}$ \textit{can not be represent as follows }%
$\omega_{n-k}=\left(  c_{1}[C_{i_{0}}\cap...\cap C_{i_{k}}]|_{C_{i_{0}}%
\cap...\cap C_{i_{k}}}\right)  \wedge\omega_{1}$\textit{\ on }$C_{i_{0}}%
\cap...\cap C_{i_{k}},$ \textit{where }$\omega_{1}\in H^{n-k-2}(C_{i_{0}}%
\cap...\cap C_{i_{k}},\mathbb{Q})$.

\textbf{Proof\ of\ Proposition }\ref{G101}: Suppose that $\omega_{n-k}\in
H^{n-k}(C_{i_{0}}\cap...\cap C_{i_{k}},\mathbb{Q})$,

\begin{center}
$\int_{\gamma}\omega_{n-k}\neq0.$
\end{center}

and $\omega_{n-k}=\left(  c_{1}[C_{i_{0}}\cap...\cap C_{i_{k}}]|_{C_{i_{0}%
}\cap...\cap C_{i_{k}}}\right)  \wedge\omega_{1}$\textit{\ on }$C_{i_{1}}%
\cap...\cap C_{i_{k}}.$ Let $\eta$ be anon zero section of the line bundle
$\mathcal{O}([C_{i_{0}}\cap...\cap C_{i_{k}}])$ on $C_{i_{1}}\cap...\cap
C_{i_{k}}$ such that the zero set of $\eta$ is exactly $C_{i_{0}}\cap...\cap
C_{i_{k}}.$ Let us consider the form

\begin{center}
$(d\log(\eta))\wedge\left(  c_{1}[C_{i_{0}}\cap...\cap C_{i_{k}}]|_{C_{i_{0}%
}\cap...\cap C_{i_{k}}}\right)  \wedge\omega_{1}$
\end{center}

on $C_{i_{1}}\cap...\cap C_{i_{k}}\backslash C_{i_{0}}\cap...\cap C_{i_{k}}.$
Let us consider the cycle $p_{1}^{-1}(\gamma)\in H_{n-k+1}(C_{i_{1}}%
\cap...\cap C_{i_{k}}\backslash C_{i_{0}}\cap...\cap C_{i_{k}},\mathbb{Q})$ as
defined in Definition \ref{t10}. Let us compute

\begin{center}
$\int_{p_{1}^{-1}(\gamma)}(d\log(\eta))\wedge\left(  c_{1}[C_{i_{0}}%
\cap...\cap C_{i_{k}}]|_{C_{i_{0}}\cap...\cap C_{i_{k}}}\right)  \wedge
\omega_{1}.$
\end{center}

Since locally around a point $w\in C_{i_{0}}\cap...\cap C_{i_{k}}$ the divisor
$C_{i_{0}}\cap...\cap C_{i_{k}}$ in $C_{i_{1}}\cap...\cap C_{i_{k}} $ is given
by $z=0,$ where $w\in U\subset C_{i_{1}}\cap...\cap C_{i_{k}},$ we see that

\begin{center}
$d\log(\eta)|_{U}=\frac{dz}{z}.$
\end{center}

From this local expression of $d\log(\eta)$ and the definition of $p_{1}%
^{-1}(\gamma)$ we deduce that

\begin{center}
$\int_{p_{1}^{-1}(\gamma)}(d\log(\eta))\wedge\left(  c_{1}[C_{i_{0}}%
\cap...\cap C_{i_{k}}]\right)  \wedge\widetilde{\omega}_{1}=2\pi\int_{\gamma
}\omega_{n-k}$
\end{center}

where $\omega_{n-k}=\left(  c_{1}[C_{i_{0}}\cap...\cap C_{i_{k}}]|_{C_{i_{0}%
}\cap...\cap C_{i_{k}}}\right)  \wedge\omega_{1}$ and $\widetilde{\omega}_{1}$
is a closed form in then tubular neighborhood of $C_{i_{0}}\cap...\cap
C_{i_{k}}$ in $C_{i_{i}}\cap...\cap C_{i_{k}}$ such that \ the restriction of
$\widetilde{\omega}_{1}$ on $C_{i_{0}}\cap...\cap C_{i_{k}}$ is $\omega_{1}.$
On the other hand since the restriction of the line bundle $\mathcal{O(}%
[C_{i_{0}}\cap...\cap C_{i_{k}}])$ on $C_{i_{1}}\cap...\cap C_{i_{k}%
}\backslash C_{i_{0}}\cap...\cap C_{i_{k}}$ is the trivial line bundle, we
deduce that the form $c_{1}[C_{i_{0}}\cap...\cap C_{i_{k}}]$ will be an exact
form on $C_{i_{1}}\cap...\cap C_{i_{k}}\backslash C_{i_{0}}\cap...\cap
C_{i_{k}}$, i.e. $c_{1}[C_{i_{0}}\cap...\cap C_{i_{k}}]=d\beta.$ This implies that

\begin{center}
$(d\log(\eta))\wedge\left(  c_{1}[C_{i_{0}}\cap...\cap C_{i_{k}}]\right)
\wedge\widetilde{\omega}_{1}=(d\log(\eta))\wedge d\beta\wedge\widetilde
{\omega}_{1}=d(d\log(\eta))\wedge\beta\wedge\widetilde{\omega}_{1}).$
\end{center}

So Stoke's Theorem implies that

\begin{center}
$\int_{p_{1}^{-1}(\gamma)}(d\log(\eta))\wedge\left(  c_{1}[C_{i_{0}}%
\cap...\cap C_{i_{k}}]\right)  \wedge\widetilde{\omega}_{1}=\int_{p_{1}%
^{-1}(\gamma)}d\left(  (d\log(\eta))\wedge\beta\wedge\widetilde{\omega}%
_{1}\right)  =\int_{\partial\left(  p_{1}^{-1}(\gamma)\right)  }(d\log
(\eta))\wedge\beta\wedge\widetilde{\omega}_{1}=0$
\end{center}

The last equality follows from the fact that $\partial\left(  p_{1}%
^{-1}(\gamma)\right)  =\emptyset.$ So we can conclude that

\begin{center}
$\int_{p_{1}^{-1}(\gamma)}(d\log(\eta))\wedge\left(  c_{1}[C_{i_{0}}%
\cap...\cap C_{i_{k}}]\right)  \wedge\widetilde{\omega}_{1}=2\pi\int_{\gamma
}\omega_{n-k}=0.$
\end{center}

On the other hand we know that

\begin{center}
$\int_{\gamma}\omega_{n-k}\neq0.$
\end{center}

So we got a contradiction. Proposition \ref{G101} is proved. $\blacksquare.$

In order to finish the proof of our Theorem we will need some facts from the
Theory of Mixed Hodge Structures.

It is easy to see that $\gamma\notin\operatorname{Im}G_{k}$ implies that there
exists $\omega_{n-k}\in H^{n-k}(C_{i_{0}}\cap...\cap C_{i_{k}},\mathbb{Q}$)
such that

\begin{center}
$\int_{\gamma}\omega_{n-k}\neq0.$
\end{center}

We may assume that $\omega_{n-k}$ is the Poincare dual of $\gamma$. The
condition $\gamma\notin\operatorname{Im}G_{k}$ implies that $\omega
_{n-k}\notin\operatorname{Im}G_{k}^{\ast}.$ So from here, Theorem \ref{D} and
Proposition \ref{G101} we conclude that we can find a form $\widetilde{\omega
}$ on $X_{0}\backslash$Sing($X_{0}$) such that $R^{[k]}(\widetilde{\omega
})=\omega_{n-k}.$ From the theory of Leray residues it follows that

\begin{center}
$\int_{p_{k}^{-1}(\gamma)}\widetilde{\omega}=\int_{\gamma}\omega_{n-k}\neq0.$
\end{center}

Lemma \ref{G1} is proved. $\blacksquare.$

\textbf{Proof\ of\ Lemma \ref{G11}:}We need to prove that if $\gamma
\notin\operatorname{Im}G_{k}$, then $\gamma_{t}=\mathcal{C}_{t}^{-1}(\gamma)$
represent a non zero class of cohomology in $H_{n}(X_{t},\mathbb{Q}).$ Let
$\omega_{t}$=$\mathcal{C}_{t}^{\ast}(\widetilde{\omega}).$ It is easy to see
from the definition of the Clemenc map that $\omega_{t}$ is a well defined
closed $n- $form on $X_{t}.$ On the other hand since $\mathcal{C}_{t}$ is
a diffeomorphism between $X_{t}\backslash\mathcal{C}_{t}^{-1}($Sing($X_{0}$))
and $X_{0}\backslash$Sing($X_{0}$) we deduce that

\begin{center}
$\int_{p_{k}^{-1}(\gamma)}\widetilde{\omega}=\int_{\gamma_{\tau}}\omega
_{t}=\int_{\gamma}\omega_{n-k}\neq0.$
\end{center}

The last inequality implies Lemma \ref{G11}. Lemma \ref{G11} is proved.
$\blacksquare.$

\textbf{Proof\ of\ Theorem\ \ref{betty1} i:}Theorem \ref{betty1} follows
directly from Theorem \ref{cl1} and Lemma \ref{G11}. Theorem \ref{betty1} part
\textbf{i} is proved. $\blacksquare.$

\textbf{Proof\ of\ Theorem\ }\ref{betty1} \textbf{ii:}Let $\alpha
_{1},..,\alpha_{k}$ be a basis of cycles of $H_{n}(\Pi(X_{0}),\mathbb{Q}),$
where $\Pi(X_{0})$ is defined in Definition \ref{POL}. From the definition of
$\Pi(X_{0})$ we can assume that the cycle $\alpha_{k}$ consists of the
$n$-dimensional simplices $S_{1_{k}},...,S_{N_{k}}$ such that the boundary of
the cycle $\alpha_{k}$ is zero. Each n dimensional simplex $S_{i}$ corresponds
to a point $q_{i}=C_{j_{0,k}}\cap...\cap C_{j_{n,k}}$ according to Definition
\ref{POL}$.$ The fact that the n-dimensional simplexes $S_{1_{k}}%
,...,S_{N_{k}}$ form a cycle means say \ that any singular points $q_{i}$ and
$q_{j}$ can be joint by Riemann surface, which means that they lie on some
$C_{j_{1,k}}\cap...\cap C_{j_{n,k}}.$ This follows directly from the fact that
the boundary of the cycle formed from $S_{1_{k}},...,S_{N_{k}}$ is zero. So on
each Riemann surface of the form $C_{j_{1,k}}\cap...\cap C_{j_{n,k}}$ that
contains the points $q_{i}$ and $q_{j}$ we can find a meromorphic form of the
third kind $\omega_{ij}$ which has poles only at the points $q_{i}$ and
$q_{j}$ with residues say $+1$ at $q_{i}$ and $-1$ at $q_{j}.$ From the
spectral sequence defined in Definition \ref{ss}, Theorem \ref{D0} and Theorem
\ref{D} we deduce that there is a holomorphic form $\omega_{k}$ in
$X_{0}\backslash$Sing($X_{0})$ such that Poincare residue of this form on each
Riemann surface of the form $C_{j_{1,k}}\cap...\cap C_{j_{n,k}}$ that contains
the points $q_{i}$ and $q_{j}$ is equal to $\omega_{ij}.$ Let the meromorphic
form $\omega_{k}$ will be non-zero on the component $C_{j_{0,k}}$ of
Sing($X_{0}),$ where $q_{i}=C_{j_{0,k}}\cap...\cap C_{j_{n,k}}$ and $q_{i}$
was defined as above. Suppose that the divisors $C_{j_{0,i}}\cap C_{j_{m,i}}$
are given by the equation $z_{m}=0$ in $C_{j_{0,i}}.$ Let us consider the
cycle $p_{n}^{-1}(q_{i})$ defined by $|z_{1}|=\varepsilon,...,|z_{n}%
|=\varepsilon$ in $C_{j_{0}}\backslash(Sing(X_{0})\cap C_{j_{0}}).$ It is easy
to see that the form $\omega_{k}$ locally around the point $q_{i}=C_{j_{0,k}%
}\cap...\cap C_{j_{n,k}}$ will be given by

\begin{center}
$\omega_{k}|_{U}=\frac{dz_{1}\wedge...\wedge dz_{n}}{z_{1}\times...\times
z_{n}}.$
\end{center}

So we can conclude that

\begin{center}
$\int_{p_{n}^{-1}(q_{i})}$ $\omega_{k}=(2\pi\sqrt[2]{-1})^{n}\neq0.$
\end{center}

From Lemma \ref{an1} we deduce that the cycles $\mathcal{C}_{t}^{-1}%
(p_{n}^{-1}(q_{i}))$ and $\mathcal{C}_{t}^{-1}(q_{i})$ are homological to each
other in $X_{t}.$ From the fact that all direct images $R^{i}\pi_{\ast}%
\omega_{\mathcal{X}/\mathcal{D}}(\log(X_{0}))$ of the sheaf 
$\omega_{\mathcal{X}/\mathcal{D}}%
(\log(X_{0}))$ are locally free sheaves on $D,$ we deduce from the exact sequence

\begin{center}
$0\rightarrow$ $\omega_{\mathcal{X}/\mathcal{D}}(\log(X_{0}))\overset{\otimes
t}{\rightarrow}\omega_{\mathcal{X}/\mathcal{D}}(\log(X_{0}))\rightarrow\omega
_{\mathcal{X}/\mathcal{D}}(\log(X_{0}))|_{X_{0}}\rightarrow0$
\end{center}

that we have the following exact sequence:

\begin{center}
$0\rightarrow H^{0}(\mathcal{X},\omega_{\mathcal{X}/\mathcal{D}}(\log(X_{0}
)))$
$\overset{\otimes t}{\rightarrow}H^{0}(\mathcal{X},\omega_{\mathcal{X}%
/\mathcal{D}}(\log(X_{0})))\rightarrow H^{0}(X_{0},\omega_{\mathcal{X}/\mathcal{D}}(\log
(X_{0}))|_{X_{0}})\rightarrow0$
\end{center}

It is easy to see that the theory of mixed Hodge structures implies that as a
free module over $H^{0}(\mathcal{D},\mathcal{O}_{\mathcal{D}})$ the module
$H^{0}(\mathcal{X},$ $\omega_{\mathcal{X}/\mathcal{D}}(\log(X_{0})))$ is of rank bigger
or equal to $\dim_{\mathbb{Q}}H_{n}(\Pi(X_{0}),\mathbb{Q)}$. So we can find
$\omega\in H^{0}(\mathcal{X},$ $\omega_{\mathcal{X}/\mathcal{D}}(\log(X_{0})))$ such that

\begin{center}
$\omega|_{D_{j_{0},k}}=\omega_{k}.$
\end{center}

Let us define $\omega_{t}=\omega|_{X_{t}}.$ So we get a holomorphic family of
holomorphic n-forms $\omega_{t,k}$ on $X_{t}$ such that

\begin{center}
$\underset{t\rightarrow0}{\lim}\omega_{t}=\omega_{0}$ and $\omega
_{0}|_{D_{j_{0,k}}}=\omega_{k}.$
\end{center}

From here we deduce that for small enouph $t$ and Lemma \ref{an1} we have that

\begin{center}
$\int_{\mathcal{C}_{t}^{-1}(q_{i})}\omega_{t}\neq0.$
\end{center}

From the last inequality we derive that the cycle $\mathcal{C}_{t}^{-1}%
(q_{i})$ is a non zero element of $H_{n}(X_{t},\mathbb{Q}).$ Now Theorem
\ref{betty1} part \textbf{ii }follows directly from Corollary \ref{betty2}.
Theorem \ref{betty1} part \textbf{ii} is proved. $\blacksquare.$

\begin{definition}
Let us define the geometric genus of $X_{t}$ as follows:
\end{definition}

\begin{center}
$p_{g}(X_{t})=\dim_{\mathbb{C}}H^{0}(X_{t},\Omega_{t}^{n})$ for $t\neq0$
\end{center}

\begin{corollary}
\label{betty3}We have the following formula for the geometric genus$:$
\end{corollary}

\begin{center}
$p_{g}(X_{t})=\sum_{i=0}^{m}p_{g}(C_{i})+\sum_{k=0}^{n-1}\sum_{|I|=k+1}
\dim_{\mathbb{C}}H^{0}(C_{I},\Omega
^{n-k})+\dim_{\mathbb{C}}H_{n}(\Pi(X_{0}),\mathbb{Q}).$
\end{center}

\section{Applications of Clemens' Theory to Complete Intersections}

\subsection{A Simple Criteria for the Existence of Jordan Block of Maximal Dimension}

In this paragraph we will prove the following Theorem, which we will apply
later to complete intersections in toric varieties:

\begin{theorem}
\label{main} Let $\pi:\mathcal{X}\rightarrow\mathcal{D}$ be a
family of K\"{a}hler varieties over the disk such that:
\end{theorem}

\begin{enumerate}
\item \textit{For t}$\neq0,$ $\pi^{-1}(t):=X_{t}$ \textit{is a non singular
variety of complex dimension }$n\geq1.$

\item $\dim_{\mathbb{C}}H^{0}(X_{t},\Omega^{n})\geq1$ \textit{for} $t\neq0. $

\item $\pi^{-1}(0)=X_{0}=\cup_{i=0}^{m}C_{i}$ \textit{is a divisor of normal
crossing and each C}$_{i}$ is irreducible.

\item $C_{0}\cap..\cap C_{n}=q$\textit{\ is a point}$.$

\item \textit{There exists a holomorphic section }$\omega\in H^{0}%
(\mathcal{X},\Omega_{\mathcal{X}/\mathcal{D}}^{n}<\log X_{0}>)$ \textit{such that
the restriction of }$\omega$ to \textit{C}$_{0},$ $\omega_{0}:=$
$\omega\left|  _{C_{0}}\right.  $ \textit{has the following expression in an
open set $\mathcal{U}$ around the singular point} $q:=C_{0}\cap..\cap C_{n}:$
\end{enumerate}

\begin{center}
$\omega_{0}\left|  _{C_{0}\cap\mathcal{U}}\right.  :=$ $\omega\left|
_{C_{0}\cap\mathcal{U}}\right.  $ $=\frac{dz_{1}\wedge...\wedge dz_{n}}%
{z_{1}\times...\times z_{n}}$
\end{center}

\textit{Then the monodromy operator of the family }$\pi:\mathcal{X}
\rightarrow\mathcal{D}$\textit{\ has a Jordan block of size
}$n+1.$

\textbf{Proof\ of\ Theorem \ref{main}:} Let $z_0,...,z_n$ be local coordinates
near $q$ such that $z_i=0$ 
local equation of for $C_{i}$, for $i=1,..,n$. Then $z_1,..,z_n$ can
be regarded as coordinates on $C_0$ near $q$. Let
$T_{0}$ be the $n$ dimensional real torus in $C_{0}\backslash\left(  C_{i}\cap
C_{0}\right)$ defined by $|z_{i}|=\varepsilon$ for $i=1,..,n.$

In the notations of Definition \ref{t1}, $T_{0}$ is just $p_n(\gamma)^{-1}(q)$,
where $\gamma$ is the homology class of the point $q$.
From the properties of the Clemens map $\mathcal{C}_{t}$ it follows
Lemma \ref{an1} that the real $n$ cycle $\mathcal{C}_{t}^{-1}(T_{0})=T_{t}$ is
homological to $\gamma_{t}:=\mathcal{C}_{t}^{-1}(\gamma).$ Let $\omega_{t}%
=\omega|_{X_{t}}.$ The idea of the proof is to show that for $t$ small enough,
$T_t$ is not homological to zero, by showing that

\begin{center}
$\int_{T_{t}}\omega_{t}\neq0.$
\end{center}

This fact together with Corollary \ref{betty2} implies Theorem \ref{main}. So
we need to prove the following Proposition:

\begin{proposition}
\label{bong1}The real cycle T$_{t}$ is not homological to zero on $X_{t}.$
\end{proposition}

\textbf{Proof\ of\ Proposition \ref{bong1}:} Condition \textbf{5} implies that
$\underset{t\rightarrow0}{\lim}\omega_{t}=\omega_{0}$ exists and locally
around the singular point $q\in C_{0}$ of $X_{0},$ $\omega_{0}$ is given by
the following expression:

\begin{center}
$\omega_{0}\left|  _{C_{0}\cap\mathcal{U}}\right.  =\frac{dz_{1}%
\wedge...\wedge dz_{n}}{z_{1}\times...\times z_{n}}.$
\end{center}

From this expression we obtain that

\begin{center}
$\int_{T_{0}}\omega_{0}=\int_{T_{0}}\frac{dz_{1}\wedge...\wedge dz_{n}}%
{z_{1}\times...\times z_{n}}=\left(  2\pi\varepsilon\sqrt{-1}\right)  ^{n}\neq0.$
\end{center}

From the fact that

\begin{center}
$\underset{t\rightarrow0}{\lim}\int_{T_{t}}\omega_{t}=\int_{T_{0}}\omega
_{0}\neq0$
\end{center}

we obtain that small $t$ we have:

\begin{center}
$\int_{T_{t}}\omega_{t}\neq0$
\end{center}

This proves Proposition \ref{bong1}. $\blacksquare.$ Theorem \ref{main} is
proved. $\blacksquare.$

\subsection{Families of Complete Intersection CY Manifolds in Projective
Space with Maximal Unipotent Monodromy}

\begin{definition}
\label{b1}We will consider the following one parameter family of complete
intersections of CY manifolds in $\mathbb{P}^{n+k}$ for $n\geq4$ and $k\geq1$
defined by the following equations$:$
\end{definition}

\begin{center}
$G_{1,t}=tF_{1}-\prod_{i=0}^{n_{1}}x_{i}=0,..,G_{k,t}=tF_{k}-\prod
_{j=n_{1}+..+n_{k-1}}^{n_{k}}x_{j}=0,$
\end{center}

\textit{where the system }$F_{1}=..=F_{k}=0$ \textit{defines a non singular CY manifolds,}

\begin{center}
\ \ \ \ \ \ \ \ \ \ \ \ \ \ $n_{i}=\deg F_{i}\geq2$ and $\sum n_{i}=n+k+1$
\end{center}

\textit{and x}$_{i}$ \textit{are the standard homogeneous coordinates in
}$\mathbb{P}^{n+k}.$

The condition $\sum n_{i}=n+k+1$ implies that the fibers $\pi^{-1}(t)=X_{t}$
for $t\neq0$ are CY manifolds of complex dimension $n$.

Let us denote the family defined in Definition \ref{b1} by

\begin{center}
$\mathcal{X\rightarrow D},$ where $\mathcal{X\subset}\mathbb{P}^{N}%
\times\mathcal{D}.$
\end{center}

\begin{theorem}
\label{bty}The family $\pi:\mathcal{X\rightarrow D}$\ has a monodromy\ T of
maximal unipotent index, i.e. (T$^{M}-id)^{n+1}=0$ and (T$^{M}-id)^{n}\neq0. $
\end{theorem}

\textbf{Proof\ of\ Theorem \ref{bty}:} Our proof is based on checking the
conditions of Theorem \ref{main} for this particular family. First we will
check that the singular fiber $X_{0}:=\pi^{-1}(0)$
consists of linear subspaces $C_{i}$ isomorphic to $\mathbb{P}^{n}.$ We will
prove that there exists a component say $C_{0}$ such that another $n$ linear
spaces $C_{i}$ intersect $C_{0}$ transversely along $\mathbb{P}^{n-1}.$ First
we will prove the following Proposition:

\begin{proposition}
\label{bty1}The singular fiber of the family $\mathcal{X\rightarrow D}$
consists of exactly $\prod_{i=1}^{k}(n_{i}-1)$ n dimensional linear subspaces
isomorphic to $\mathbb{P}^{n}.$
\end{proposition}

\textbf{Proof\ of\ Proposition \ref{bty1}: }From the Definition \ref{b1} it
follows that $X_{0}:=\pi^{-1}(0)$ is defined by the equations:

\begin{center}
$\prod_{i=0}^{n_{1}}x_{i}=0,..,\prod_{j=n_{1}+..+n_{k-1}}^{n_{k}}x_{j}=0$
\end{center}

and thus $X_{0}$ is a union of linear subspaces of dimension $n$ defined by
the equations:

\begin{center}
$x_{j_{1}}=0,..,x_{j_{k}}=0,$
\end{center}

where $0\leq x_{j_{1}}\leq n_{1},..,n_{1}+..+n_{k-1}<x_{j_{k}}\leq N.$ From
these equations Proposition \ref{bty1} follows directly. $\blacksquare.$

\begin{proposition}
\label{bty2}There exists a component $C_{0}$ of the singular fiber $X_{0}$ of
the family $\pi:\mathcal{X\rightarrow D}$ such that $C_{0}$ is isomorphic to
$\mathbb{P}^{n}.$ Let ($z_{0}:..:z_{n})$ be the homogeneous coordinates in
$\mathbb{P}^{n},$ then there exists $n$ components of the singular fiber
$X_{0}$ say $C_{1},..,C_{n}$ such that $C_{i}\cap C_{0}$ is defined in
$D_{0}\approxeq\mathbb{P}^{n}$ by the equation $z_{i}=0.$
\end{proposition}

\textbf{Proof\ of\ Proposition \ref{bty2}: }Let us denote by

\begin{center}
$d_{1}=n_{1},$ $d_{2}:=n_{1}+n_{2},..,d_{j}:=n_{1}+...+n_{j},..,d_{k}%
=n_{1}+...+n_{k}.$
\end{center}

We already proved that each component of $X_{0}$ is a linear subspace
isomorphic to $\mathbb{P}^{n}.$ Let $D_{0}$ be the subspace defined by the equations:

\begin{center}
$x_{0}=x_{d_{1}}=x_{d_{2}}=...=x_{d_{k-1}}=0.$
\end{center}

The coordinates $(x_{i_{1}}:x_{i_{2}}:...:x_{n+k})$ define a homogeneous
coordinates on $C_{0}$ if all $x_{i_{m}}$ are different from the $x_{d_{j}}$
that define $C_{0}.$ Let us define $C_{i_{m}}$ for $1\leq m\leq n$ by the
following system of equations: $x_{j_{1}}=..=x_{j_{k}}=0,$ where all the
indexes $(j_{1},..,j_{k})$ with an exception of one is equal to indexes
$(0,d_{1},..,d_{k-1}).$ These equations define a component of the singular
fiber. Indeed the necessary and sufficient condition for a system of linea
equations $x_{j_{1}}=..=x_{j_{k}}=0$ to define a component of the singular
fiber $X_{0}$ of the family $\pi:\mathcal{X\rightarrow D}$ is $0\leq j_{1}\leq
d_{1},$ $d_{1}<j_{2}\leq d_{3},..,d_{k-1}<j_{k}\leq d_{k}.$ Here we are using
the conditions that $n_{i}\geq2$ and $n_{1}+...+n_{k}=n+k.$ Clearly $C_{i_{m}%
}\cap C_{0}$ will be a $n-1$ dimensional linear subspace in $C_{0}%
\approxeq\mathbb{P}_{n-1}.$ From here it follows that we constructed $C_{i}$
for $i=1,..,n$ such that $C_{0}\cap...\cap C_{n}=(0:...:1)=q.$ So we proved
Proposition \ref{bty2}. $\blacksquare.$

\begin{proposition}
\label{bty3}There exists a section $\omega\in H^{0}(\mathcal{X},\Omega
_{\mathcal{X}/\mathcal{D}}^{n}\log<X_{0}>)$ such that if $\omega|_{X_{t}}=\omega_{t}$
for $t\neq0$ then:
\end{proposition}

\begin{center}
$\underset{t\rightarrow0}{\lim}\omega_{t}=\omega_{0}$
\end{center}

\textit{where when we restrict }$\omega_{0}$ \textit{on }$C_{0}\approxeq
\mathbb{P}^{n},$ \textit{then }$\omega_{0}\left|  _{U_{i}}\right.  $\textit{is
given by the formula:}

\begin{center}
$\omega_{0}|_{U_{i}}=\frac{dz^{1}\wedge..\wedge dz^{n}}{z^{1}...z^{n}}.$
\end{center}

\textit{where }$U_{i}$ \textit{is the standard covering of }$\mathbb{P}^{n}.$

\textbf{Proof\ of\ Proposition \ref{bty3}:} Since we choose $p=(0,0,..,1)$ we
will work in the open set $U_{n+k}:=\{(x_{0}:..:x_{n+k})|x_{n+k}\neq0\}.$ Let
us consider the meromorphic form

\begin{center}
$\Omega_{t}:=\frac{dz^{1}\wedge..\wedge dz^{n+k}}{\tilde{G}_{t,1}..\tilde
{G}_{t,k}},$ where $\tilde{G}_{t,i}:=G_{t,i}(\frac{x_{0}}{x_{n+k}}%
,..,\frac{x_{n+k-1}}{x_{n+k}},1)$ and $z^{i}=\frac{x_{i}}{x_{n+k}}.$
\end{center}

on $U_{n+k}\times\mathcal{D}.$ By taking $k$ times the Leray residue of
$\Omega_{t}$ we define the holomorphic n-form $\omega_{t}$ on $X_{t}.$ Suppose
that \ $C_{0}$ in $U_{n+k}\times\mathcal{D}$ is given by the equations (by
suitably reordering the coordinates):

\begin{center}
$z^{n+1}=..=z^{n+k}=0$ and $t=0.$
\end{center}

(See \cite{L}.) Let $z^{1},..,z^{n}$ be the rest of the coordinates in
$U_{n+k}.$ Let us choose an open set $\mathcal{W}$ in $\mathcal{X\subset
}\mathbb{P}^{N}\times D^{\ast}$ where $z^{1},..,z^{n},t$ are local coordinates
of $\mathcal{W}$, then \ from the definition of the Leray residue it follows
that on $\mathcal{W}$ we have:

\begin{center}
$\omega_{t}|_{\mathcal{W}}:=\frac{dz^{1}\wedge..\wedge dz^{n}}{\partial
_{n+1}...\partial_{n+k}(\tilde{G}_{t,1}\times..\times\tilde{G}_{t,N-n})},$
where $\partial_{i}=\frac{\partial}{\partial z^{i}}.$
\end{center}

From the last formulas we get that

\begin{center}
$\underset{t\rightarrow0}{\lim}\partial_{n+1}...\partial_{n+k}(\tilde{G}%
_{t,1}\times..\times\tilde{G}_{t,k})=z^{1}\times..\times z^{n}.$
\end{center}

So we deduce that

\begin{center}
$\underset{t\rightarrow0}{\lim}\omega_{t}|_{\mathcal{W}}=\omega_{0}%
=\frac{dz^{1}..\wedge dz^{n}}{z^{1}\times..\times z^{n}}.$
\end{center}

Proposition \ref{bty3} is proved. $\blacksquare.$

Propositions \ref{bty2} and \ref{bty3} imply that the conditions of Theorem
\ref{main} are satisfied by the family $\pi:\mathcal{X\rightarrow D}.$ So
Theorem \ref{main} implies Theorem \ref{bty}. Theorem \ref{bty} is proved.
$\blacksquare.$

\begin{remark}
\label{yau}It is not at all difficult to generalized the contsruction of the
one parameter family of complete intersections of CY manifold in
$\mathbb{CP}^{N}$ with maximal unipotent element in the monodromy group to the
case of of complete intersections of general type. More precisely it is very
easy to prove that there exists a family $\pi:\mathcal{Y}\rightarrow
\mathcal{D}$ of complete intersections in $\mathbb{CP}^{N}$ of general type
such that if the dimension of the fiber is $n,$ then the monodromy operator
$T$ contains exactly $p_{g}(Y_{t})$ Jordan blocks of dimension $n+1,$ where
$p_{g}=\dim_{\mathbb{C}}H^{0}(Y_{t},\Omega_{Y_{t}}^{n}).$
\end{remark}

\section{Construction of a Family of CY Manifolds in Toric Varieties with
Maximal Unipotent Monodromy}

\subsection{Introduction to Toric Varieties}

Let $T=(\mathbb{C}^{\ast})^{N}.$ Let $\mathbb{N}$ be a rank N lattice, and
$\Sigma$ be a complete fan relative to $\mathbb{N}$. (See \cite{O}.)

Notations:

\begin{description}
\item $\mathbb{P}_{\Sigma}=V:$ toric variety defined by $\Sigma.$

\item $\mathbb{N}^\vee:$ lattice dual to $\mathbb{N}.$

\item $\Sigma(k):$ the set of k-cones in $\Sigma.$

\item $\tau^{\vee}:$ the dual of a cone $\tau.$

\item $D_{\rho}:$ the toric divisor corresponding to $\rho\in\Sigma(1).$

\item $\mathcal{O}(D)=\mathcal{O}_V(D):$ the line bundle (invertible
sheaf) associated to the divisor $D$.
\end{description}

We follow the usual combinatorial description of a cone $\tau$ in $\Sigma,$
and often use the set of edges of $\tau$ or its primitive generators to denote
the cone. Let $O_{\rho}$, $\rho\in\Sigma(1),$ be the codimension 1 $T$-orbit
and $D_{\rho}$ be their respective closures. Thus they are irreducible
$T$-invariant Weil divisors in $V.$ Let $R$ be the polynomial ring
$\mathbb{C}[x_{\rho},$ $\rho\in\Sigma(1)].$ If we declare that $\deg x_{\rho
}=[D_{\rho}]\in A_{N-1}(V),$ then $R$ becomes a\ $A_{N-1}(V)$-graded ring,
where $A_{N-1}(V)$ is the Picard group of $V.$ We denote the degree $[D]$
subspace in $R$ by $R_{[D]}.$ It is known that 
(see \cite{Au}\cite{Cox} 
\cite{Mu}):

\begin{center}
$R_{[D]}\approxeq H^{0}(V,\mathcal{O}_{V}(D)).$
\end{center}

The graded ring $R$ is called the homogeneous coordinate ring of V. The
isomorphism above can be described as follows. The dual lattice $\mathbb{N}%
^{\vee}$ can be viewed as the lattice of characters of the group $T.$
We denote a character multiplicatively corresponding to $\nu\in\mathbb{N}%
^{\vee}$ by $\chi^{\nu}.$  Then

\begin{center}
$H^{0}(V,\mathcal{O}_{V}(D))=
\oplus_{\nu\in P_{D}\cap\mathbb{N}%
^{\vee}}\mathbb{C\cdot\chi}^{\nu}$
\end{center}

where

\begin{center}
$P_{D}=\{\nu\in\mathbb{N}^{\vee}\otimes_{\mathbb{Z}}\mathbb{R}%
|\left\langle \nu,\rho\right\rangle \geq-a_{\rho}$ $\forall\rho\}.$
\end{center}

Here $\rho$ denotes the primitive generators of the 1-cones, and $D=\sum
a_{\rho}D_{\rho}$ is the Weil divisor defining the line bundle $\mathcal{O}%
_{V}(D).$ The isomorphism $\phi_{D}$ above is given by $\phi_{D}:\chi^{\nu
}\rightarrow\prod_{\rho}x_{\rho}^{\left\langle \nu,\rho\right\rangle +1}.$ The
assertion above gives a description of the space of sections of all
equivariant line bundles over $V.$

Let us consider the section $x_{\rho}$. Note that $\nu=0$ is in $P_{D_{\rho}%
},$ hence by the isomorphism above $\phi_{D_{\rho}}(\chi^{0})=x_{\rho}$ is a
section of $\mathcal{O}_{V}(D).$ In fact the zero set of $x_{\rho}=0$ is
exactly the divisor $D_{\rho}.$

For mirror symmetry, the most interesting case is when the convex hull of the
1-cone generators $\rho\in\mathbb{N}^{\vee}\otimes_{\mathbb{Z}%
}\mathbb{R}$ form a reflexive polytope $\Delta,\,$\ and that $D$ is the
anticanonical divisor. (See \cite{Bat} and \cite{Bor}.) In this case $P_{D}$ is
the polar dual of $\Delta.$ Here are some examples:

\begin{enumerate}
\item $V=\mathbb{P}^{2},$ $D=D_{1}+D_{2}+D_{3},$ D$_{i}$ being the $i^{th}%
-$hyperplane. The homogeneous coordinate ring of $V$ is the usual
$\mathbb{C}[x_{1},x_{2},x_{3}]$. If we
identify $A_{1}(V)$ with $\mathbb{Z}$ such that $[H]\rightarrow1$, then
$x_{i}$ has degree 1. The polytope $\Delta$ is generated by the 1-cone
generators $\rho_{i}$ is the triangle with vertices (1,0), (0,1), (-1,-1). Its
dual $\nabla$ is the triangle with vertices (2,-1), (-1,2), (-1,-1). It is
easy to check that $\nabla$ has exactly 10 lattice points $m.$ They are in 1-1
correspondence with the degree 3 monomials in $x_{i}.$

\item Let $\xi$ be a primitive third root of unity.
Let $\mathbb{Z}_3$ acts on $\mathbb{P}^2$ by $[x_1,x_2,x_3]\mapsto
[\xi x_1, \xi^{-1} x_2, x_3]$. Resolve all the
singularities in $\mathbb{P}^2/\mathbb{Z}_3$ "minimally".
The result is a toric variety $V$ with
9 hyperplanes $D_{i}$ all together (with two linear relations). The
anticanonical divisor in $V$ is $D=D_{1}+...+D_{9}.$ The polytope $\Delta$
generated by the 1-cone generators $n_{\rho}$ is now the triangle with
vertices (2,-1), (-1,2), (-1,-1). Its dual $\nabla$ is the triangle with
vertices (1,0), (0,1), (-1,-1). Since $\nabla$ has exactly 4 lattice points,
$K_{V}^{-1}$ has 4 sections -- one of them being $x_{1}\cdots x_{9}.$
\end{enumerate}

\subsection{Construction of One Parameter Family of CY Complete Intersection
in a Toric Variety with Maximal Jordan Block in the Monodromy Operator}

Let $\pi_{1},...,\pi_{k}$ be a partition of the set $\Sigma(1).$ We assume
that $F_{i}\in H^{0}(V,\mathcal{O}_{V}(\sum_{\rho\in\pi_{i}}D_{\rho}))$
together define a nonsingular subvariety $X$ of codimension k:

\begin{center}
$F_{1}=...=F_{k}=0.$
\end{center}

Put $n=N-k.$ Since the anticanonical class is

\begin{center}
$[K_{V}^{-1}]=\underset{\rho}{\sum}D_{\rho},$
\end{center}

by adjunction $X$ is a Calabi-Yau manifold. We consider the following
1-parameter family of complete intersections $X_{t}$ defined by

\begin{center}
$G_{i,t}=tF_{i}-\underset{\rho\in\pi_{i}}{\prod}x_{\rho}=0,$ $i=1,...,k.$
\end{center}

Let us denote this family when $|t|<1$ by $\pi:\mathcal{X}\rightarrow
\mathcal{D}.$ The fiber $X_{0}$ is the union of toric subvarieties

$$C_{\sigma}:=\cap_{\rho\in\sigma}D_{\rho}$$

where $\sigma$ is any subset of $\Sigma(1)$ such that
$\vert$%
$\sigma\cap\pi_{i}|=1$ for all $i.$ It turns out that (see Appendix)
$C_{\sigma}$ is nonempty iff the 1-cones in $\sigma$ generate a cone in
$\Sigma.$ In this case, $C_{\sigma}$ is a nonsingular irreducible toric
subvariety in $V$ of codimension $k$ corresponding to the k-cone $\sigma
\in\Sigma(k).$ Thus such a nonempty $C_{\sigma}$ is a component in $X_{0}.$

\begin{remark}
The assumption that the $F_{i}$ defines a codimension k subvariety is
important. This put a strong constraint on the kind of partitions $\pi
_{1},...,\pi_{k}$ that are allowed. A general partition of $\Sigma(1)$ will
fail to satisfy this condition. Here is an example. Take the 4-dimensional
weight projective space $\mathbb{P}[9,6,1,1,1].$ After a minimal
desingularization, the resulting toric variety has 9 1-cones $\rho
_{1},...,\rho_{9}$ in its fan $\Sigma.$ After suitable ordering, it is easy to
find a partition of the form $\pi_{1}=\{\rho_{1},\rho_{2}\},$ $\pi_{2}%
=\{\rho_{3},\rho_{4}\},$ $\pi_{3}=\{\rho_{5},...,\rho_{9}\}$ such that
$\{\rho_{1},\rho_{3}\},$ $\{\rho_{1},\rho_{4}\},$ $\{\rho_{2},\rho_{3}\},$
$\{\rho_{2},\rho_{4}\},$ are all primitive collections. In this case, $X_{0}$
would be empty because the $C_{\sigma}=D_{\rho_{i}}\cap D_{\rho_{j}}\cap
D_{\rho_{k}},$ for $\rho_{i}\in\pi_{1},$ $\rho_{j}\in\pi_{2},$ $\rho_{k}\in
\pi_{3},$ will all be empty.
\end{remark}

\textit{Throughout this section, we shall make the
following additional assumption: that
the convex hull of the set of primitive
generators }$\rho$ \textit{of the 1-cones in }$\Sigma$\textit{\ \ is a
reflexive polytope, which we shall denote by
}$\Delta.$

\begin{theorem}
\label{BY}The monodromy operator of the family $\mathcal{X}\rightarrow
\mathcal{D}$ \ has one Jordan block of dimension $n+1.$
\end{theorem}

\subsection{Proof of Theorem \ref{BY}}

\begin{lemma}
\label{BY1}Let $C_{\sigma_{0}}$ be a fixed component in $X_{0}.$ Then there
exists n other components $C_{\sigma_{1}},...,C_{\sigma_{n}}$ with the
property that $C_{\sigma_{0}}\cap C_{\sigma_{i}}$ is a codimension $k+1$ toric
subvariety, and that
\end{lemma}

\begin{center}
$C_{\sigma_{0}}\cap...\cap C_{\sigma_{n}}$
\end{center}

\textit{is a toric fixed point.}

\textbf{Proof:} Since the fan $\Sigma$ is complete, we can find an
$N-$ cone, say $\tau$, contaning $\sigma_{0}$ as a k-face. Write
$\sigma_{0}=\left\langle \rho_{1},...,\rho_{k}\right\rangle ,$ where the
$\rho_{i}\in\pi_{i}$ are the canonical generators of $\sigma_{0}.$ Similarly,
write $\tau$ in terms of its 1-cone generators:

\begin{center}
$\tau=\left\langle \rho_{1},...,\rho_{k},\rho_{1}^{^{\prime}},...,\rho
_{n}^{^{\prime}}\right\rangle .$
\end{center}

Then each $\rho_{j}^{^{\prime}}$ lies in a unique say $\pi_{i_{j}}.$ Since $V$
is non-singular, $\tau$ is a simplicial cone, so that any k of its generators
generate a k-cone in $\Sigma.$ In particular we have the k-cones

\begin{center}
$\sigma_{j}:=\left\langle \rho_{1},...,\hat{\rho}_{i_{j}},...,\rho_{k}%
,\rho_{j}^{^{\prime}}\right\rangle ,$ $j=1,...,n.$
\end{center}

That is, $\sigma_{j}$ is obtained from $\sigma_{0}$ by replacing the generator
$\rho_{i_{j}}$ of $\sigma_{0}$ by $\rho_{j}^{^{\prime}}.$ Since both
$\rho_{i_{j}}$ and $\ \rho_{j}^{^{\prime}}$ live in the same $\pi_{i_{j}},$ it
follows that $C_{\sigma_{j}},$ as defined above, is a component in $X_{0}.$ By construction,

\begin{center}
$C_{\sigma_{0}}\cap C_{\sigma_{j}}=\cap_{\rho\in\sigma_{0}\cup\sigma_{j}%
}D_{\rho}.$
\end{center}

But $\sigma_{0}\cup\sigma_{j}$ is the list $\sigma_{0}$ adjoint with $\rho
_{j}^{^{\prime}},$ hence gives a ($k+1)$-cone in $\Sigma.$ So $C_{\sigma_{0}%
}\cap C_{\sigma_{j}}$ is a codimension $(k+1)$ toric subvariety in $V.$ Moreover,

\begin{center}
$C_{\tau}=\cap_{\rho\in\sigma_{0}\cup...\cup\sigma_{n}}D_{\rho}=\cap_{\rho
\in\tau}D_{\rho}$
\end{center}

which is codimension $k+n=N$ toric subvariety of $V.$ Hence it is a fixed
point.$\blacksquare.$

Note that the nonsingular toric subvariety $C_{\sigma_{0}}$ comes with a
standard affine coordinates on the affine patch containing fixed point above.
Namely, they are obtained from restricting the standard affine coordinates on
the patch

\begin{center}
$U_{\tau}=Hom_{sg}(\tau^{\vee},\mathbb{C}^{\ast})\approxeq(\mathbb{C}^{\ast
})^{N}.$
\end{center}

Here the isomorphism is determined by the choice of ordering of the set of
primitive generators of the cone $\tau^{\vee}.$

\begin{lemma}
\label{BY2}There exists a meromorphic form on $V=\mathbb{P}_{\Sigma}$ with
simple poles along each hypersurface $F_{i}=0.$
\end{lemma}

\textbf{Proof:} First let $X$ be a complex $N$-fold, and $K_{X}$ its canonical
bundle. Let $\{U_{\sigma}\}$ be a covering of charts on $X,$ whose coordinates
$U_{\sigma}\rightarrow\mathbb{C}^{N}$ we denote $z^{\sigma}=(z_{1}^{\sigma
},...,z_{N}^{\sigma}).$ From this data, we get a frame $dz^{\sigma}%
=dz_{1}^{\sigma}\wedge...\wedge z_{N}^{\sigma}$ on each $U_{\sigma}$ for the
bundle $K_{X},$ hence a dual frame for the dual bundle $K_{X}^{-1}.$ Suppose
$f$ is a nonzero global section of $K_{X}^{-1}.$ Then relative to a dual frame
above, $f$ is represented as a holomorphic function $f_{\sigma}(z^{\sigma}).$
Note that on any overlap $U_{\sigma}\cap U_{\tau}\neq\emptyset,$ the functions
$f_{\sigma},$ $f_{\tau}$ transform by the same transition function as the
frame $dz^{\sigma}$ and $dz^{\tau}.$ It follows that the local expressions

\begin{center}
$\frac{dz^{\sigma}}{f_{\sigma}(z^{\sigma})}$
\end{center}

together define a global meromorphic $N$-form on $X$ with poles along the zero
locus $f=0.$ We now assume that $X$ is an algebraic variety and the
$U_{\sigma}$ are affine patches. 
We now apply this to an $N$-dimensional complete
toric variety $V=\mathbb{P}_{\Sigma},$ as before. 

Recall that $\mathbb{P}_{\Sigma}$ is covered by the affine subvarieties
$U_{\sigma}\approxeq\mathbb{C}^{N},$ labelled by $N$-cones $\sigma$ in
$\Sigma.$ Here the isomorphism is determined by an ordering of the set of $N$
integral generators of the cone $\sigma.$ We denote by $z^{\sigma}%
=(z_{1}^{\sigma},...,z_{N}^{\sigma})$ the coordinates of this isomorphism. Let

\begin{center}
$\nabla:=\{\nu\in\mathbb{N}^{\vee}\otimes\mathbb{R}|\left\langle
\nu,x\right\rangle \geq-1,$ $\forall x\in\Delta\}.$
\end{center}

This is the dual of the reflexive polytope $\Delta.$

The global sections of the anticanonical bundle $K_{V}^{-1}$ on $V=\mathbb{P}%
_{\Sigma}$ has the following description. There is a basis of $H^{0}%
(V,K_{V}^{-1})$ which corresponds 1-1 with the set $\nabla\cap\mathbb{N}%
^{\vee}.$ Let's fix an ordering of the set $\Sigma(1)$ and denote the
generators by $\rho_{1},...,\rho_{S}\in\mathbb{N}.$ Let $L$ be the kernel of
the natural map

\begin{center}
$\mathbb{Z}^{S}\rightarrow\mathbb{N},$ $m\rightarrow\sum m_{i}\rho_{i}.$
\end{center}
Let $L^\perp$ be the orthogonal complement of $L$ in $\mathbb{Z}^S$
with respect to the standard inner product $l\cdot l'$. 
Then $\mathbb{N}\approxeq\mathbb{Z}^{S}/L,$ and the natural pairing $L^{\perp
}\times\mathbb{Z}^{S}/L\rightarrow\mathbb{Z},$ $(l,l^{^{\prime}}+L)\mapsto
l\cdot l^{^{\prime}},$ defines a canonical
isomorphism $L^{\perp}\approxeq\mathbb{N}^{\vee}.$  Then
a general section of $K_{V}^{-1}$ is of the form

\begin{center}
$f=x_{1}...x_{S}\underset{\nu\in\nabla\cap\mathbb{N}^{\vee}}{\sum
}c_{\nu}x^{\nu}.$ 
\end{center}

(See \cite{Au}, \cite{Cox} and \cite{Mu}.)
Here the $c_{\nu}$ are arbitrary complex numbers, $x_{i}$ are the homogeneous
coordinates of $V,$ viewed as a section of $\mathcal{O}(D_{\rho_{i}}),$ and
$x^{\nu}$ is the monomial in the $x_{i}$ with exponents $\nu\in L^{\perp
}\subset\mathbb{Z}^{S}.$

As an example, when $V=\mathbb{P}^{N}$, we have $K_{V}^{-1}=\mathcal{O}%
_{\mathbb{P}^{N}}(N+1),$ and $x_{i}$ are the usual homogeneous coordinates. In
this case a general section of $K_{V}^{-1}$ above is exactly a degree $N+1$
homogeneous polynomial in the $x_{i}.$

We now return to the general case. Let $e_{1},...,e_{N}$ denotes the standard
basis of $\mathbb{Z}^{N}$. Let $\sigma$ be any $N$-cone in the fan $\Sigma,$
and fix an ordering $\nu_{1},...,\nu_{N}$ for the primitive generators of the
dual cone $\sigma^{\vee}.$

\begin{proposition}
\label{BY3}Relative to the frame dual to $dz^{\sigma}$ above, sections of
$K_{V}^{-1}$ are represented by polynomial functions of the form
\end{proposition}

\begin{center}
$f_{\sigma}(z)=z_{1}...z_{N}\underset{v\in\nabla_{\sigma}}{\sum}c_{v}z^{v}.$
\end{center}

\textit{Here }$z_{i}=z_{i}^{\sigma},$ $c_{v}$\textit{\ are arbitrary numbers,
}$\nabla_{\sigma}\subset\mathbb{Z}^{N}$\textit{\ is the image of }$\nabla
\cap\mathbb{N}^{\vee}$\textit{\ under the isomorphism }$\mathbb{N}%
^{\vee}\rightarrow\mathbb{Z}^{N}$\textit{\ \ determined by }$\nu
_{i}\rightarrow e_{i}$\textit{.}

\textbf{Proof:}A priori\textbf{, }$f_{\sigma}(z)$ given above is a Laurent
poly nomial in $z.$ We will first show that it is in fact a polynomial. Recall
that every element $\nu\in\nabla$ satisfies $\left\langle \nu,x\right\rangle
\geq-1$ for all $x\in\Delta.$ Let $v$ be the image of $\nu$ under the
isomorphism $\mathbb{N}^{\vee}\rightarrow\mathbb{Z}^{N}$ determined by
$\nu_{i}\rightarrow e_{i},$ and let $a$ be the preimage of $x\in\Delta$ under
the dual map $\mathbb{Z}^{N}\rightarrow\mathbb{N}$. Then $v\cdot
a=\left\langle \nu,x\right\rangle \geq-1.$ In particular if $x$ is then
preimage generator of the cone $\sigma,$ then $\left\langle \nu_{k}%
,x\right\rangle =1$ for one $k$ and zero otherwise. This means that the
preimage $a$ of $x$ is such that $e_{i}\cdot a=1$ if $i=k$ and zero otherwise,
i.e. $a=e_{k}.$ Clearly, every $e_{k}$ can be realized as the preimage of some
primitive generator $x\in\Delta$ of $\sigma.$ This shows that for any
$v\in\nabla_{\sigma},$ we have $v\cdot e_{k}=\left\langle \nu,x\right\rangle
\geq-1$ for all $k.$ It follows that every Laurent monomial $z_{1}%
...z_{N}\cdot z^{v}$ appearing in $f_{\sigma}(z)$ has a nonnegative exponent,
i.e. $f_{\sigma}(z)$ is a
polynomial function on $\mathbb{C}^{N}.$ 

Let $\tau$ be any other $N$-cone in the fan $\Sigma,$ and fix an ordering
$\mu_{1},...,\mu_{N}$ for the primitive generators of $\tau^{\vee}.$
Clearly the respective generators $\nu_{i},\mu_{j}$ for $\sigma^{\vee%
},\tau^{\vee},$ determine an $A\in(a_{ij})\in\mathbb{GL}(N,\mathbb{Z})
$ such that

\begin{center}
$\nu_{i}=\underset{j}{\sum}a_{ij}\mu_{j}.$
\end{center}

For simplicity, we write $z=z^{\sigma},$ $w=z^{\tau}.$ The coordinates
transition function is then given by

\begin{center}
$z_{i}=w^{A_{i}},$ $A_{i}=(a_{i1},...,a_{iN}).$
\end{center}

We consider the region of overlap of $U_{\sigma}\cap U_{\tau}$ where
$z_{i}\neq0$ and $w_{j}\neq0$ for all $i,j.$ First we show that the frames
$dz,$ $dw$ of $N$-forms are related by the transformation law

\begin{center}
$dz=\det(A)\frac{z_{1}...z_{N}}{w_{1}...w_{n}}dw.$
\end{center}

We have

\begin{center}
$\frac{dz}{z_{1}...z_{N}}=d\log z_{1}\wedge...d\log z_{N}.$
\end{center}

Now $\log z_{i}=\sum_{j}a_{ij}\log w_{j}$ on the overlap. So the
transformation law above follows immediately.

It remains to show that there is a polynomial function $f_{\tau}(w)$ such that

\begin{center}
$f_{\sigma}(z)=\det(A)\frac{z_{1}...z_{N}}{w_{1}...w_{N}}f_{\tau}(w).$
\end{center}

The matrix $A$ defines an isomorphism $B:\mathbb{R}^{N}\rightarrow\mathbb{R}%
^{N},$ $Bv=\sum_{j}a_{ij}v_{i}.$ Then

\begin{center}
$z^{v}=w^{Bv}$ .
\end{center}

Moreover the image of $\nabla_{\sigma}\subset\mathbb{R}^{N}$ under $B$ is
$\nabla_{\tau}.$ It follows that

\begin{center}
$f_{\sigma}(z)=z_{1}...z_{N}\cdot\underset{v\in\nabla_{\sigma}}{\sum}%
c_{v}z^{v}=z_{1}...z_{N}\underset{u\in\nabla_{\tau}}{\sum}c_{B^{-1}u}w^{u}.$
\end{center}

So if we set

\begin{center}
$f_{\tau}(w):=w_{1}...w_{N}\cdot\underset{u\in\nabla_{\tau}}{\sum}%
\det(A)^{-1}c_{B^{-1}u}w^{u},$
\end{center}

then we obtain the desired function. $\blacksquare.$

Proposition \ref{BY3} implies Lemma \ref{BY2}. $\blacksquare.$

As an example, when $V=\mathbb{P}^{N},$ we can choose the fan $\Sigma$ in
$\mathbb{Z}^{N},$ so that the 1-cone primitive generators are $e_{0}%
=-e_{1}-\cdots-e_{N}, e_1,..,e_N.$ In this case, we have $N+1$ standard affine patches
$U_{\sigma_{i}},$ where $\sigma_{i}$ is the $N$-cone generated by
$e_{0},..,\hat{e}_{i},...,e_{N}.$ The affine coordinate are $z^{\sigma_{i}%
}=(\frac{x_{0}}{x_{i}},...,1_{i},...\frac{x_{N}}{x_{i}}).$ For simplicity, we
label each cone $\sigma_{i}$ simply by the integer $i.$ Then the coordinate
transition functions are given by

\begin{center}
$z_{i}^{\sigma}=z_{j}^{\tau}\frac{x_{\tau}}{x_{\sigma}},$ $i\neq\sigma,\tau;$
$z_{\tau}^{\sigma}=\frac{1}{z_{\sigma}^{\tau}}.$
\end{center}

The transition functions for coordinate $N$-forms are given by

\begin{center}
$dz^{\sigma}=(-1)^{\sigma+\tau}\frac{x_{\tau}}{x_{\sigma}}dz^{\tau}.$
\end{center}

We now prove the general toric version of Proposition \ref{bty3}.

\begin{proposition}
\label{BY4}There exists a family of holomorphic n-forms $\omega_{t}$ on
$X_{t}$ for $t\neq0$ such that
\end{proposition}

\begin{center}
$\underset{t\rightarrow0}{\lim}\omega_{t}=\omega_{0}$
\end{center}

\textit{where }$\omega_{0}$\textit{\ is a section on }$X_{0}$\textit{\ with
the following properties: restricted to the component }$C_{\sigma_{0}}%
$\textit{, }$\omega_{0}$\textit{\ has the form:}

\begin{center}
$\omega_{0}|_{U_{\varepsilon}}=\frac{dz_{1}\wedge...dz_{n}}{z_{1}...z_{n}} $
\end{center}

\textit{where (}$z_{1},...,z_{n})$\textit{\ \ are coordinates on the affine
patch }$U_{\varepsilon}\subset C_{\sigma_{0}}.$

\textbf{Proof:} Let $U_{\tau}\approxeq\mathbb{C}^{N}$ be the affine patch in
$V=\mathbb{P}_{\Sigma}$ as given in the proof of Proposition \ref{BY3}%
.\thinspace Let $z^{\tau}=z=(z_{1},...,z_{N})$ be the affine coordinates. Here
we choose the coordinates so that the affine subvariety $U_{\varepsilon
}=U_{\tau}\cap C_{\sigma_{0}}$ is defined by

\begin{center}
$z_{n+1}=...=z_{N}=0$
\end{center}

and the $z_{1},...,z_{n}$ are affine coordinates of $U_{\varepsilon}.$ Choose
a small analytic neighborhood
$W\subset\mathcal{X}\cap U_{\tau}\times \mathcal{D}$ so that
$z_{1},...,z_{n},t$ are coordinates on $X_{t}\cap W$ for small $t.$

Clearly

\begin{center}
$f:=G_{t,1}...G_{t,k}\in H^{0}(V,\mathcal{O}(\underset{\rho\in\Sigma(1)}{\sum
}D_{\rho}))=H^{0}(V,K_{V}^{-1}).$
\end{center}

By Lemma \ref{BY2}, we have a meromorphic $N$-form given by

\begin{center}
$\Omega_{t}=\frac{dz}{g_{t}(z)}$
\end{center}

where

\begin{center}
$g_{t}(z)=z_{1}...z_{N}(1-tf_{1}(z))...(1-tf_{k}(z)).$
\end{center}

Here $f_{i}(z)$ is a Laurent polynomial representing the section $F_{i}$ of
$\mathcal{O}(\sum_{\rho\in\pi_{i}}D_{\rho})$ on the algebraic torus
$z_{1}...z_{N}\neq0.$ Note that for all $t,$ $g_{t}$ is polynomial in the $z$
by construction. Integrating the form $\Omega_{t}$ via Leray residue $k$
times, the result is a section $\omega_{t}$ of the sheaf of $n$-forms, whose
restriction on $X_{t}\cap W$ is given by

\begin{center}
$\omega_{t}=\frac{dz_{1}\wedge...\wedge dz_{n}}{\frac{\partial}{\partial
z_{n+1}}...\frac{\partial}{\partial z_{N}}g_{t}(z)}$ .
\end{center}

It is clear that the dominator of $\omega_{t}$ goes to $z_{1}...z_{n}$ as
$t\rightarrow0.$ This proved our assertion. $\blacksquare.$

Theorem \ref{BY} now follows directly from Theorem \ref{main}. $\blacksquare.$

\begin{remark}
\label{yau2}One can generalize the same constructions to the case of compete
intersections of general type in toric varieties for which the canonical class
is very ample. Namely one can prove that there exists one parameter family of
complete intersections of general type $\mathcal{Y}\rightarrow\mathcal{D}$ for
which the generic fibers is of dimension $n$ and the canonical class is very
ample, such that the monodromy operator contains exactly $p_{g}(Y_{t})$ Jordan
blocks of dimension $n+1.\,$
\end{remark}

\section{Comparisons}

\subsection{Relationship with Hypergeometric Functions}

In \cite{HLY}, the construction of the point of maximal unipotent monodromy
has been done for the family of CY hypersurfaces consisting of all
anticanonical hypersurfaces modulo equivalence under $T$ in a fixed toric
variety. The construction there uses an entirely different approach. Namely,
one exploits the fact that the periods of the family are solutions to a system
of linear partial differential equations, known as a GKZ hypergeometric system
and its generalization. (See \cite{HLY} and the references there.) The
construction is then done by carefully analyzing the structure of certain
singularities of the PDE system. Explicit computation of certain period is
also necessary in that approach.

We now briefly compare the approach of \cite{HLY} with the present approach.
For simplicity, we consider the case of hypersurfaces in an $N$ dimensional
toric variety $V=\mathbb{P}_{\Sigma}.$ Recall that in order to apply Theorem
\ref{main}, we must construct a meromorphic $N$-form $\Omega_{t}$ on $V$ with
a simple pole along the CY hypersurface $X_{t}:$

\begin{center}
$G_{t}:=tF-\underset{\rho}{\prod}x_{\rho}=0.$
\end{center}

Here $F\in H^{0}(V,K_{V}^{-1})$ is a fixed section. In the affine coordinates
$z^{\sigma}=(z_{1},...,z_{N})$ on an affine patch $U_{\sigma}\approxeq
\mathbb{C}^{N},$ $\Omega_{t}$ takes the form:

\begin{center}
$\Omega_{t}=\frac{dz_{1}\wedge...\wedge dz_{N}}{z_{1}...z_{N}(1-tf(z))}$
\end{center}

where $f(z)$ is a Laurent polynomial representing the section $F.$ The Leray
residue of $\Omega_{t}$ gives a holomorphic $n$-form on $X_{t}.$ We then
integrate $\omega_{t}$ over a real $n$-torus, giving a holomorphic function of
$t$ bounded near $t=0.$ This function is a period of the CY manifold $X_{t}$
near $t=0.$ It is easy to see that this period can also be obtained by
integrating $\Omega_{t}$ over a real $N$-torus $\alpha_{0},$ $|z_{1}%
|=...=|z_{N}|=\varepsilon,$ in $T:$

\begin{center}
$\int_{\alpha_{0}}\Omega_{t}.$
\end{center}

But this is precisely a special value of the holomorphic solution of GKZ
system given in \cite{HLY}. In fact, in \cite{HLY}, $G_{t}$ is replaced by a
general section, and the special Laurent polynomial $1-tf(z)$ is replaced by
its general counterpart. The result is a multi-variable holomorphic function
(the variables being the coefficients of the general Laurent polynomial). A
power series expansion of this function can be easily computed. (See
\cite{HLY}.) The period $\int_{\alpha_{0}}\Omega_{t}$ above can be obtained by
specializing the multi-variable function to 1-variable function by restricting
it to the one-parameter family $G_t=0$.

\subsection{Comparison of Maximal Unipotent Monodromy
and the SYZ Conjecture}

Clemens' theorem \ref{cl1} and our construction
of the maximal unipotent monodromy operator $T$ for
complete intersection CY manifolds in a toric variety
show that there exist
special cycles, 
namely the invariant cycle $\gamma$ and the cycles $\alpha
_{1},...,\alpha_{n},$ such that 

\begin{center}
$T(\alpha_{n})=\gamma+\alpha_{1}+...+\alpha_{n}.$
\end{center}

From Clemens' theory, we know that the cycle $\gamma$ 
can be realized as a $n$ dimensional torus. 
The cycle $\alpha_{n}$ can probably be realized
as an $n$ dimensional sphere.
A theorem of Clemens says that $\gamma$ and $\alpha_n$
have intersection number 1.

On the other hand, the SYZ conjecture asserts that a CY $n$-fold
should admit the structure of a fibration over an $n$-sphere
with generic fibers given by special Lagrangian $n$-tori. Moreover this fibration
admits a canonical section. Thus the section
is an $n$-cycle having intersection number 1 with 
the fiber cycle. It is quite clear that the cycle $\gamma$ and $\alpha_n$ above
should be realized precisely as the $n$-sphere section and the $n$-torus fiber
respectively. 
This general relation deserves further investigation.
For families of polarized K3 surfaces having maximal unipotent monodromy, 
this can be verified by using
isometric deformation (relative to the Calabi-Yau metric) 
of the surface, as we will discuss next.

\def\func{}

\begin{theorem}
\label{JT}
Let $\pi:\mathcal{X}\rightarrow\mathcal{D}$ be a family of polarized
K3 surfaces with polarization class $e\in H^{2}(X,\mathbb{Z})$ such
that $<e,e>=2n$ for a fixed positive integer $n.$ Suppose that $\pi
^{-1}(0)=X_{0}$ is a singular surface. Fix $t\neq0$
and suppose the monodromy operator $T$ acting on
the second homology group of $X=\pi^{-1}(t)$ is such that
$(T-id)^{3}=0$ and $(T-id)^{2}\neq0.$ 
Let $g_{t}$ be the
Calabi-Yau metrics on $X=X_t$ such that $[Im~g_{t}]=e$. 
Let $\gamma,$ $\alpha_{1}$ and $\alpha_{2}$ be cycles such that
$$T(\gamma)=\gamma,~~T(\alpha_{1})=\gamma+\alpha_{1},~~T(\alpha
_{2})=\gamma+\alpha_{1}+\alpha_{2},$$
$$<e,\gamma>=<\gamma,\gamma>=<\alpha_{2},\alpha_{2}>=0,~~<\gamma,\alpha_{2}>=1.$$
Then there exists a torus fibration 
$\psi:X\rightarrow S^{2}$ such that the generic fiber
$\psi^{-1}(s)$ is a Lagrangian 2-torus
with respect to $Re~g_{t}$ representing the homology class $\gamma$.
Moreover, there is a section $\sigma:S^2\rightarrow X$ such that $\sigma(S^2)\subset X$
is a Lagrangian submanifold
representing the homology class $\alpha_2-\gamma$.
\end{theorem}

Before we prove the theorem, let's construct a family
of polarized K3 surface that fulfills the
conditions stated in Theorem \ref{JT}. Let $\mathcal{E}_{t}:y^{2}=x(x-1)(x-t)$
be a family of elliptic curves when $t\in\mathcal{D},$ where $\mathcal{D}%
:=\{t\in\mathbb{C}|$ $|t|<1\}$ and $\mathcal{E}_{t}\subset\mathbb{CP}%
^{2}\times\mathcal{D}.$ Then $\mathcal{E}_{t}\times_{\mathcal{D}}%
\mathcal{E}_{t}\rightarrow\mathcal{D}$ is a family of 
abelian surfaces over the unit disc. Let $\pi:\mathcal{K}_{t}\rightarrow
\mathcal{D}$ be the family of Kummer surfaces associated with this family of
abelian surfaces. In \cite{JT} it was proved that the family $\pi
:\mathcal{K}_{t}\rightarrow\mathcal{D}$ satisfies the conditions of Theorem
\ref{JT}. See page 252-253 in \ref{JT}. On page 252 an explicit construction
of the cycles $\gamma$ and $\alpha_{2}$ is given and moreover from Clemens'
theory it follows that

\begin{center}
$<\gamma,\gamma>=<\alpha_{2},\alpha_{2}>=0$ and $<\gamma,\alpha_{2}>=1$. 
\end{center}

In particular, $\gamma$ and $\alpha_{2}$ span a hyperbolic lattice

\begin{center}
$\mathbb{H}=\left(
\begin{array}
[c]{ll}%
0 & 1\\
1 & 0
\end{array}
\right)  .$
\end{center}

It is a well known fact that for a K3 surface $X$ 
$H_{2}(X,\mathbb{Z})$ is a lattice isomorphic to $\Lambda_{K3}:=\mathbb{H}^{3}\oplus
(-E_{8})^{2}.$ Let us define (transcendental lattice)

\begin{center}
$T_{e}:=\{v\in\Lambda_{K3}|<e,v>=0\}.$
\end{center}

We know also that $T_{e}\approxeq\mathbb{Z} e^{\ast}
\oplus\mathbb{H}^{2}\mathbb{\oplus%
}(-E_{8})^{2},$ where $<e^{\ast},e^{\ast}>=-<e,e>.$ So we may suppose that $%
\ T_{e}=\mathbb{Z}e^{\ast}\oplus\{\gamma,\alpha _{2}\}
\mathbb{\oplus H\oplus}%
(-E_{8})^{2}.$ See \cite{JT}. In \cite{JT} it is shown that the cycles $%
\gamma$ and $\alpha_{2}\in T_{e}.$

\textbf{PROOF of Theorem \ref{JT}:} 
We will use isometric deformation with respect to the CY metric $g_{t}$
on a K3 surface whose class of cohomology $[\func{Im}g_{t}]=e$ and the
epimorphism of the period map for K3 surfaces to prove that K3 surafce can
be realized as a Lagrangian fibration with respect to each CY metric $g_{t}.$
The Lagrangian fibration that we will construct will have the properties
that the homology class of the generic fibre will be $\gamma $ and the
homology class of the base will be $\delta =\alpha _{2}-\gamma .$ For the
description and applications of isometric deformations see \cite{To}.
Theorem \ref{JT} will follow directly from the following Lemma:

\begin{lemma}
\label{JT2} \textbf{i. }Let $Y$ be a K3 surface whose Neron-Severi group $%
NS(Y)$ spanned by a cycle $\gamma $ and such that $<\gamma ,\gamma >=0,$
then $\gamma $ can be realized as a non singular elliptic curve $C$ on $Y$
and the linear system $|C|$ defines an elliptic fibration $|C|:Y\rightarrow 
\mathbb{CP}^{1}.$ \textbf{ii } Let $Y$ be a K3 surface whose Neron Severi group 
$NS(Y)$ spanned by a cycle $\delta $ and such that $<\delta ,\delta >=-2,$
then either $\delta $ or $-\delta $ can be realized as an embedded complex
line $\mathbb{CP}^{1}\subset Y.$
\end{lemma}

\textbf{PROOF: }First we will prove part \textbf{i} of Lemma \ref{JT2}. From
the condition that Neron Severi group $NS(Y)$ is spanned by cycles $\gamma $
such that $<\gamma ,\gamma >=0$, we deduce that there exists a line bundle $%
\mathcal{L}_{\gamma },$ whose Chern class $c_{1}(\mathcal{L}_{\gamma
})=\gamma .\,$\ From Riemann-Roch theorem we deduce that

\begin{center}
$\dim H^{0}(Y,\mathcal{L}_{\gamma})+\dim H^{2}(Y,\mathcal{L}_{\gamma})\geq2.$
\end{center}

Serre's duality implies that $\dim H^{2}(Y,\mathcal{L}_{\gamma })=\dim
H^{2}(Y,\mathcal{L}_{\gamma }^{\ast }),$ where $\mathcal{L}_{\gamma }^{\ast
} $ is the dual line bundle of $\mathcal{L}_{\gamma }.$ From here we obtain
that either $\dim H^{0}(Y,\mathcal{L}_{\gamma })>1$ or $\dim H^{0}(Y,%
\mathcal{L}_{\gamma }^{\ast })>1.$ So either $\gamma $ or $-\gamma $ can be
realized as an effective divisor on $Y.$ 
Recall

\begin{theorem} (\cite{Sh} p559.)
\label{Sh}If an effective divisor $D>0$ on a K3 surface Y satisfies the
conditions $<D,D>=0$ and $<D,E>\geq 0$ for any effective divisor $E>0,$ then
the linear system $|D|$ contains a divisor of the form $mC,$ where $m>0$ and 
$C$ is an elliptic curve.
\end{theorem}
Since $NS(Y)$\ is spanned by a 
cycle $\gamma $ such that $<\gamma ,\gamma >=0$,
it follows that
that $\gamma $ or $-\gamma $ can be realized a non sigular elliptic
curve $C.$ From the standard exact sequence on $C$:

\begin{center}
$0\rightarrow \mathcal{O}_{Y}(-C)\rightarrow \mathcal{O}_{Y}\rightarrow 
\mathcal{O}_{C}\rightarrow 0$,
\end{center}
we have, in the associated long exact sequence,

\begin{center}
$H^{0}(\mathcal{O}_{Y})\rightarrow H^{0}(\mathcal{O}%
_{C})\rightarrow H^{1}(\mathcal{O}_{Y}(-C))\rightarrow H^{1}(\mathcal{O}%
_{Y})=0.$
\end{center}
Since the restriction map $H^{0}(\mathcal{O}_{Y})=\mathbb{C}%
\rightarrow H^{0}(\mathcal{O}_{C})=\mathbb{C}$ is an isomorphism, we 
have $H^{1}(\mathcal{O}_{Y}(-C))=0.$ From Serre's duality applied to
K3 surafces we get $H^{1}(\mathcal{O}_{Y}(-C))=H^{1}(\mathcal{O}%
_{Y}(C))=0$ and $H^{2}(\mathcal{O}_{Y}(C))=0.$ It is a standard fact that
the linear system $|C|$ defines a map $\phi :Y\rightarrow \mathbb{CP}^{1}.$ For
details see page 560 of \cite{Sh}. Part \textbf{i} of Lemma \ref{JT2}
is proved.

Part \textbf{ii} of Lemma \ref{JT2} follows directly from Corollary 3 on
page 560 of \cite{Sh}. Lemma \ref{JT2} is proved. $\blacksquare .$

One of the main consequences of the existence of a Calabi Yau metric on a K3
surface X with a holomorphic form $\omega_{X}$ is the fact that the
covariant derivative of the holomorphic form $\omega_{X}$ is zero, i.e. $%
\nabla \omega_{X}=0.$ This implies that the following three forms $\{\func{Re%
}\omega_{X},\func{Im}\omega_{X}$ and $\func{Im}(g)\}$ are parallel, closed
and non degenerate forms on $X$. It is easy to see that:

\begin{center}
$\left\langle \func{Re}\omega_{X},\func{Im}\omega _{X}\right\rangle
=\left\langle \func{Re}\omega_{X},\func{Im}g\right\rangle =\left\langle 
\func{Im}g,\func{Im}\omega _{X}\right\rangle =0.$
\end{center}

Here $\func{Im}(g)$ is the imaginary part of the Calabi Yau metric on $X.$
If $\left\| \func{Re}\omega _{X}\right\| ^{2}=\left\| \func{Im}\omega
_{X}\right\| ^{2}=\left\| \func{Im}(g)\right\| ^{2}=1,$ then these three
forms define three integrable complex structures $I,J$ and $K$ on $X$ such
that $I^{2}=J^{2}=K^{2}=-id$ and $IJ+JI=IK+KI=JK+KJ=0.$ For all these facts
see \cite{To}. The isometric deformation is defined as a new integrable
complex structure $aI+bJ+cK$ on $X,$ where $a^{2}+b^{2}+c^{2}=1.$
We'll denote
this new complex structure on $X$ by $Y$. In \cite{To} it was
proved that Re$\omega _{Y}=A(\func{Re}\omega _{X}),$ $\func{Im}\omega _{Y}=A(%
\func{Im}\omega _{X}),$ where $A\in SO(3).$ The class $A(\func{Im}g)$ will
be a class of cohomology of type $(1,1)$ on $Y$ and will define a new Calabi
Yau metric on $Y$ which is isometric to the Calabi Yau metric on $X$ as a
Riemannian metrics. In \cite{To} we proved that there exists $A\in SO(3)$
such that

\begin{center}
$\int_{\gamma }A(\func{Re}\omega _{X})=\int_{\gamma }A(\func{Im}\omega
_{X})=0.$
\end{center}

Let $Y$ and $C$ be
defined as in Lemma \ref{JT2}. Notice that when we do an isometric
deformation, we are not changing either the $C^{\infty }$ structure 
or the
Riemannian structure on $Y.$ So we may consider that $C$ realizes the cycle $%
\gamma $ as an embedded two dimensional torus in $X.$

\begin{lemma}
\label{iso0} $C$ is a Lagrangian cycle on $X$
\end{lemma}

\textbf{PROOF:} We need to prove the following two facts: \textbf{1. }The
volume form of the restriction of the CY metric $g_{X}$ on $C$ is given by
the following expression: 
\begin{equation}
Vol(g_{X}|_{C})=a \func{Re}(\omega _{X})+b \func{Im}(\omega _{X})
\label{0}
\end{equation}
and \textbf{2.} that
\begin{equation}
\func{Im}(g_{X}|_{C})=0\mbox{ .}  \label{00}
\end{equation}
Since $C$ is an elliptic curve in $Y$, it follows that $\omega
_{Y}|_{C}=0,$ so that $\func{Re}\omega _{Y}|_{C}=\func{Im}\omega _{Y}|_{C}=0.$
The volume form
of the restriction of the Ricci flat Riemannian metric on $C$ is equal to $%
\func{Im}(g_{Y})|_{C}$, ie.
\begin{equation}
Vol(g_{X}|_C)
=\func{Im}(g_{Y})|_{C}.
\label{2}
\end{equation}
So from (\ref{2}) and the general facts about isometric deformation we deduce
that 
\begin{equation}
Vol(g_{X}|_{C})=a \func{Re}(\omega _{X})|_{C}+b \func{Im}(\omega
_{X})|_{C}+c \func{Im}(g_{X})|_{C}\mbox{ .}  \label{3}
\end{equation}
where $a,b,c$ are some real numbers 
given by
\begin{equation}
a =\int_{C}\func{Re}(\omega _{X}),\ \ b =\int_{C}\func{Im}%
(\omega _{X}),\ \ c =\int_{C}\func{Im}(g_{X})\mbox{ .}  \label{4}
\end{equation}

Since $\gamma \in T_{e}$ and the cohomology class $[Im~g_{X}]=e$ we can
conclude that: 
\begin{equation}
0=\int_{\gamma }e=\int_{C}\func{Im}(g_{X})=c.  \label{5}
\end{equation}
This proves fact \textbf{1. } We now prove
fact \textbf{2,} i.e. that the restriction of the form 
$\func{Im}(g_{X})$ on $C
$ is identically zero. This follows directly from the following formula,
proved in \cite{To}: 
\begin{equation}
\func{Im}(g_{X})|_{C}=c \func{Im}(g_{X})|_{C}\ \text{ and }\ c
=\int_{C}\func{Im}g_{X}\mbox{ .}
\end{equation}

So Lemma \ref{iso0} is proved. $\blacksquare.$

\begin{remark}
\label{iso1}Exactly in the same way we can prove that the cycle $\delta
=\alpha_{2}-\gamma$ can be realized as a Lagrangian sphere. Indeed we can
deform isometrically the complex structure on $X$ to $Z$ with respect to the
CY metric $g_{X}$ so that $\delta$ will be realized as a complex projective
line $\mathbb{P}^{1}$ embedded in $Z.$ Repeating the arguments from Lemma \ref
{iso0} we deduce that $\mathbb{P}^{1}$ is a Lagrangian sphere embedded in $X.$
\end{remark}

\textbf{The End of the Proof of Theorem \ref{JT}: }Indeed we realized the
cycles $\gamma $ and $\alpha _{2}-\gamma $ as complex analytic curves on $Y$
and on $Z.$ Then we know from Theorem \ref{Sh} that $Y$ is an elliptic
fibration whose ''generic'' fibre is an elliptic curve. Lemma \ref{iso0}
implies that the fibres of this fibration are Lagrangian submanifolds. Since
we know that $\left\langle \gamma ,\delta \right\rangle =1,$ we can conclude
that the basis of the Lagrangian fibration is the Lagrangian sphere that
realizes the cycle $\delta .$ $\blacksquare .$

\section{Appendix: Clemens' Cell Complex for Toric Hypersurfaces}

As before let $\Delta$ be the convex hull of generators $\rho$ of the
$1$-cones in $\Sigma.$ This is an $N$ dimensional polytope. But it has more.
It comes equipped with a canonical simplicial decomposition induced by
$\Sigma.$ In particular the boundary $\partial\Delta$ is topologically an
$N$-sphere which comes equipped with a simplicial decomposition.

\begin{lemma}
\label{BY5}Let $\tau$ be a set of primitive generators of 1-cones in $\Sigma.$
Put $k=|\tau|.$ The following are equivalent:
\end{lemma}

\textbf{i. }$\cap_{\rho\in\tau}D_\rho$\textit{is non empty.}

\textbf{ii. }\textit{the cone is generated by }$\tau$ \textit{is a k-cone in
}$\Sigma.$

\textbf{iii. }\textit{the convex hull of }$\tau$\textit{\ is a }%
$(k-1)$\textit{-cell in }$\partial\Delta.$

\textbf{Proof:}(\textbf{ii})$\Leftrightarrow$(\textbf{iii}) is obvious. We
will show that (\textbf{i}) and (\textbf{ii}) are equivalent. Recall that
there is an order reversing 1-1 correspondence between T-orbits $O_{\sigma}$
and cones $\sigma$ in $\Sigma$, and that the closure of V($\sigma$) of $O_{\sigma}$
is given by

\begin{center}
V($\sigma$)=$\underset{\gamma\supset\sigma}{\coprod}O_{\gamma}.$
\end{center}

Note that $D_{\rho}=V(\rho)$ for $\rho\in\Sigma(1)$. Thus

\begin{center}
$\cap_{\rho\in\tau}D_{\rho}=\underset{\gamma\supset\tau}{\coprod}O_{\gamma}.$
\end{center}

If this is nonempty, then we have some cone $\gamma\in\Sigma$ containing $\tau.$
Since the toric variety V is nonsingular, $\gamma$ is a simplicial cone. this
implies that any collection of k edges of $\gamma$ generates a k-face of
$\gamma.$ So $\tau$ generates a k-face, which we also call $\tau,$ of
$\gamma.$ In particular $\tau\in\Sigma.$ Conversely, if $\tau$ generates a
cone in $\Sigma$, obviously every $D_{\rho},$ $\rho\in\tau,$ contains
$O_{\rho}.$ In this case, $\cap_{\rho\in\tau}D_{\rho}$ is nonempty.
$\blacksquare.$

Note that if the intersection $\cap_{\rho\in\tau}D_{\rho}$ is nonempty, then
it has dimension $N-|\tau|.$ Moreover, the correspondence above between
nonempty intersections of the $D_{\rho}$'s and the cells of $\partial\Delta$
is order reversing. Namely, if two sets $\tau,$ $\tau^{\prime}$ of primitive
generators yield nonempty intersections, then $\cap_{\rho\in\tau}D_{\rho
}\subset\cap_{\rho\in\tau^{\prime}}D_{\rho}$ iff the $conv(\tau)\supset
conv(\tau^{\prime}).$

We now specialize our family $\mathcal{X}\rightarrow\mathcal{D}$ to the case
of hypersurfaces

\begin{center}
$tF-\underset{\rho}{\prod}x_{\rho}=0.$
\end{center}

Thus $N=n+1.$

\begin{theorem}
The Clemens polytope $\Pi(X_{0})$ is a simplicial complex which is naturally
isomorphic to $\partial\Delta.$
\end{theorem}

\textbf{Proof:} Recall that $X_{0}=\cup_{\rho}D_{\rho}.$ By definition of
$\Pi(X_{0}),$ each $D_{\rho}$ corresponds to a vertex in $\Pi(X_{0}).$ For
$\rho\neq\rho^{\prime},$ $D_{\rho}\cap D_{\rho^{\prime}}$ corresponds to
1-cell if $D_{\rho}\cap D_{\rho^{\prime}}$ has dimension $n-1=N-2,\,$and so
on. Thus given a set $\tau$\ of 1-cones with $|\tau|=k,\cap_{\rho\in\tau
}D_{\rho}$ corresponds to $(k-1)$-cell in $\Pi(X_{0})$ iff $\cap_{\rho\in\tau
}D_{\rho}$ has dimension $n-k+1.$ This is an order reversing correspondence
between $n-k+1$ dimensional intersections $\cap_{\rho\in\tau}D_{\rho}$ of
$D_{\rho}$'s and $(k-1)$-cells of the simplicial complex $\partial\Delta$ as
simplicial complex. $\blacksquare.$

\begin{corollary}
$H_{n}(\Pi(X_{0}),\mathbb{Q})\approxeq\mathbb{Q}.$
\end{corollary}

In the case when $\mathcal{X}$ is a family of complete intersections in
$\mathbb{P}_{\Sigma},$ one can still define
a cell complex $\Pi(X_{0})$ similar to the definition of
Clemens' complex.
It is expected that the cell complex $\Pi(X_{0})$ will again be isomorphic to
n-sphere. We have verified this for many examples. It turns out that the
complex is typically non-simplicial, unlike in the case of hypersurfaces. It
is an interesting combinatorial problem to describe the complex $\Pi(X_{0})$
in simple terms.

\end{document}